\documentclass[12pt,a4paper,oneside]{article}
\usepackage[displaymath, mathlines]{lineno}
\usepackage{authblk}
\usepackage{xr}
\usepackage{natbib}
\usepackage{amsmath}
\usepackage{mathtools}
\usepackage{mathrsfs}  
\usepackage{amsfonts}
\usepackage{amssymb}
\usepackage{amsthm}
\usepackage{graphicx}
\usepackage{subcaption}
\usepackage{soul, color, colortbl}
\usepackage{booktabs}
\usepackage[margin=1.25in]{geometry}
\usepackage{enumerate}
\usepackage[shortlabels]{enumitem}
\usepackage{url}
\usepackage{afterpage}
\usepackage{multirow}
\usepackage{xr}

\newtheorem{lemma}{Lemma}
\newtheorem{proposition}{Proposition}
\newtheorem{remark}{Remark}

\newcolumntype{H}{>{\setbox0=\hbox\bgroup}c<{\egroup}@{}}

\begin{document}

\title{ {\bf Multivariate Location and Scatter Matrix Estimation Under Cellwise and Casewise Contamination } }
\author[1]{{\bf Andy Leung}}
\author[2]{{\bf Victor J. Yohai}}
\author[1]{{\bf Ruben H. Zamar}}
\affil[1]{{\small Department of Statistics, University of British Columbia, 3182-2207 Main Mall, Vancouver, British Columbia V6T 1Z4, Canada}}
\affil[2]{{\small Departamento de Matem\'atica, Facultad de Ciencias Exactas y Naturales, Universidad de Buenos Aires, Ciudad Universitaria, Pabell\'on 1, 1426, Buenos Aires, Argentina}}
\renewcommand\Authands{, }

\maketitle

\begin{abstract}
We consider the problem of multivariate location and scatter matrix estimation when the data contain  cellwise and casewise outliers. \citet{agostinelli:2014} propose a two-step approach to deal with this problem: first apply a univariate filter to remove cellwise outliers and second   apply a generalized S-estimator to downweight casewise outliers.   We improve  this proposal in three main directions. First, we introduce a consistent bivariate filter to be used in combination with the univariate filter in the first step. Second, 
 we propose a  new  fast subsampling procedure to generate starting points for the generalized S-estimator in the second step. Third, we consider a non-monotonic weight function for the generalized S-estimator to better deal with casewise outliers in high dimension.  A simulation study and real data example show that, unlike the original two-step procedure, the modified  two-step approach performs and scales well for high dimension. Moreover, the modified procedure outperforms the original one and  other state of the art robust procedures under cellwise and casewise data contamination.
\end{abstract}

\section{Introduction}

In this paper, we address the  problem of robust estimation of  multivariate location and scatter matrix under cellwise and casewise contamination. 

Traditional robust estimators assume a casewise contamination model for the data where the majority of the cases are  assumed to be free of contamination. Any case that deviates from the model distribution is then  flagged as an outlier. In situations where only a small number of cases are contaminated this approach works well. However, if a small fraction of cells in a data table are contaminated but in such a way that a large fraction of cases are affected, then traditional robust estimators may fail. This problem, referred to as propagation of cellwise outliers, has been discussed by \citet{alqallaf:2009}.  Moreover, as pointed out by  \citet{agostinelli:2014} both types of data contamination, casewise and cellwise, may occur together.  

 Naturally, when  data  contain both cellwise and casewise outliers, the problem becomes more difficult. To address this problem, \citet{agostinelli:2014} proposed a two-step procedure: first, apply a univariate filter (UF) to  the data matrix $\mathbb X$ and set the flagged cells to missing values, NA's; and second, apply the generalized S-estimator (GSE) of \citet{danilov:2012} to the incomplete data set. Here, we call this two-step procedure UF-GSE. It was shown in \citet{agostinelli:2014} that UF-GSE is simultaneously robust against cellwise and casewise outliers. However, this procedure has three limitations, which are addressed in this paper:
\begin{itemize}
\item The univariate filter does not handle  well moderate-size cellwise outliers. 
\item The GSE procedure used in the second step loses robustness against casewise outliers for  $p > 10$.
\item  The initial estimator EMVE used in the second step does not scale well to higher dimensions  ($p > 10$).
\end{itemize} 

 \citet{rousseeuw:2015} pointed out that  to filter the  variables based solely on their value may be too limiting as no correlation with other variables is taken into account. A not-so-large contaminated cell that  passes the univariate filter could be flagged when viewed together with other correlated components, especially for highly correlated data.  To overcome this deficiency,  we introduce a consistent bivariate filter  and use it in combination with UF and a new filter developed by \citet{rousseeuw:2016} in the first step of the two-step procedure.

\citet{maronna:2015a} made a remark that UF-GSE, which uses a fixed loss function $\rho$ in the second step, cannot handle well high-dimensional casewise outliers. S-estimators with a fixed loss function exhibit an increased Gaussian efficiency when $p$ increases, but at the same time lose their robustness \citep[see ][]{rocke:1996}. Such curse of dimensionality has also been observed for UF-GSE in our simulation study.  To overcome this deficiency, we constructed a new robust estimator called {\it Generalized Rocke S-estimator} or {\it GRE} to replace GSE in the second step. 

The first step of filtering is generally fast, but the second step is slow due to the computation of the extended minimum volume ellipsoid (EMVE), used as initial estimate by the generalized S-estimator. The standard way  to compute EMVE  is by subsampling, which  requires an impractically large number of subsamples when $p$ is large, making the computation extremely slow.
To reduce the high computational cost of the two-step approach in high dimension,  we  introduce a new subsampling procedure based on clustering. The initial estimator computed in this way is  called EMVE-C.

 The rest of the paper is organized as follows.  In Section \ref{sec:GSE-bivariate-filtering}, we describe some existing filters and introduce a new consistent bivariate filter. By consistency, we mean that, when $n$ tend to infinity and the data do not contain outliers, the proportion of  data points flagged by the  filter tends to zero. We also show in Section \ref{sec:GSE-bivariate-filtering} how the bivariate filter  can be used in combination with the other filters in the first step. In Section \ref{sec:GSE-GRE}, we introduce the GRE to be used in place of GSE in the second step. In Section \ref{sec:GSE-computing-issue}, we discuss the computational issues faced by  the initial estimator, EMVE,   and introduce a new cluster-based-subsampling procedure called EMVE-C.   In Section \ref{sec:GSE-MCresults} and \ref{sec:GSE-example}, we compare the original and modified two-step approaches with several state-of-the-art  robust procedures in  an extensive simulation study. We also give there a real data example. Finally, we conclude in Section \ref{sec:GSE-conclusions}.  The Appendix contains all the proofs. We also give a separate  document called ``Supplementary Material", which contains further details,  simulation results, and other related material.

\section{Univariate and Bivariate Filters}\label{sec:GSE-bivariate-filtering}

Consider a random sample of $\mathbb X = (\pmb X_1,\dots,\pmb X_n)^t$, where  $\pmb{X}_{i}$ are first generated from a central parametric distribution, $H_0$, and then some cells, that is, some entries in  $\pmb{X}_{i}=(X_{i1},\dots,X_{ip})^{t}$ , may be independently contaminated.   A {\it filter} $\mathcal{ F}$  is a procedure that flags cells in a data table and replaces them by NA's. Let $f_n$  be the fraction of cells in the data table  flagged by the filter. A {\it consistent filter} for a given distribution $H_0$ is one that asymptotically will not flag any cell if the data come from  $H_0$. That is, $lim_{n \rightarrow \infty} f_n =0$  a.s. $[H_0]$.

\begin{remark}\label{remark:filter-definition}
  Given a collection of filters $\mathcal{ F}_1,...,\mathcal{ F}_k$ they can be combined in  several ways: (i) they can be {\it united}  to form a new filter, $ \mathcal{ F}_U =  \mathcal{F}_1\cup \cdots \cup \mathcal{ F}_k$, so that the resulting filter, $ \mathcal{ F}_U $, will flag  all the cells flagged by at least  one  of them;  (ii) they can be {\it intersected}, so that  the resulting filter, $ \mathcal{ F}_I =  \mathcal{F}_1 \cap \cdots \cap \mathcal{ F}_k$, will only  flag the cells identified by all of them; and (iii) 
a filter, $\mathcal F$, can be conditioned to yield a new filter, $\mathcal F_C$,
so that $\mathcal F_C$ will only filter the cells filtered by $\mathcal F$ which  satisfy a given condition $C$.
\end{remark}

\begin{remark}\label{remark:filter-consistency} 
 It is clear that  $ \mathcal{ F}_U$ is a consistent filter provided all the filters $ \mathcal{ F}_i$, $i=1,\dots,k$ are  consistent filters. On the other hand, $ \mathcal{ F}_I$ is a consistent filter provided at least one of the filters  $ \mathcal{ F}_i$, $i=1,\dots,k$ is a consistent filter.   Finally, it is also clear that if  $\mathcal{ F}$   is a consistent filter, so is $\mathcal{ F}_C$.
\end{remark}

We describe now three basic filters, which  will be later combined to obtain a powerful consistent filter for use  in the first step of our two-step procedure.

\subsection{A Consistent Univariate Filter (UF) }

This is the initial filter introduced in  \citet{agostinelli:2014}. Let $X_1, \dots, X_n$ be a random (univariate) sample of observations. Consider a pair of initial location and dispersion estimators, $T_{0n}$ and $S_{0n}$, such as the  median and median absolute deviation (MAD) as adopted   in this paper. Denote the standardized sample by $Z_{i} = (X_{i} -  T_{0n})/S_{0n}$. Let $F$ be a chosen reference distribution for $Z_{i}$. Here, we use the standard normal distribution, $F = \Phi$.

Let ${F}^+_{n}$ be the empirical distribution function for the absolute standardized value, that is, 
\begin{linenomath} 
\[
	{F}^+_{n}(t) = \frac{1}{n}\sum_{i=1}^n I( |Z_{i}| \le t).
\]
The proportion of flagged outliers  is defined by 
\begin{equation}\label{eq:2SGS-GY-d}
\begin{aligned}
 d_{n}  &=\sup_{t\ge \eta} \left\{ F^+
(t)-{F}^+_{n}(t) \right\}^{+}, 
\end{aligned}
\end{equation}
\end{linenomath} 
where $\{a \}^+$ represents  the positive part of $a$, $F^+$ is the distribution of $|Z|$ when $Z\sim F$, and $\eta = (F^+)^{-1}(\alpha)$ is a large quantile of $F^+$. We use $\alpha = 0.95$ for univariate filtering as the aim is to detect large outliers, but other choices could be considered. Then, we flag $\lfloor nd_{n} \rfloor$ observations with the largest absolute standardized value, $|Z_i|$, as cellwise outliers and replace them by NA's.

The following proposition states this is a consistent filter. That is,  even when the actual distribution is unknown, asymptotically, the univariate filter will  not  flag  outliers when the tail of the chosen reference distribution is heavier than (or equal to)  the tail of the actual distribution. 

\begin{proposition}[\citealp{agostinelli:2014}]\label{prop:2SGS-GY-asymptotic} 
Consider a random variable $X \sim F_{0}$ with $F_0$ continuous.
Also, consider a pair of  location and dispersion estimators $T_{0n}$ and $S_{0n}$ such that $T_{0n}\rightarrow\mu_0 \in \mathbb R$ and $S_{0n}\rightarrow\sigma_0 > 0$ a.s. [$F_0$].  Let $F_0^+(t) = P_{F_0}( |\frac{X - \mu_0}{\sigma_0}| \le t)$.
 If  the reference distribution $F^+$  satisfies the inequality
\begin{linenomath} 
\begin{equation}\label{eq:2SGS-reference}
\max_{t\ge \eta}\left\{F^+(t)-F_{0}^+(t) \right\}\leq 0,
\end{equation} 
\end{linenomath} 
then
\begin{linenomath} 
\[
\frac{n_{0}}{n}\rightarrow0\text{ a.s.,}%
\]
\end{linenomath} 
where 
\begin{linenomath} 
\[
n_{0}= \lfloor nd_{n} \rfloor.
\]
\end{linenomath} 
\end{proposition}

We define the global univariate filter, UF, as the union of all the consistent filters described above, applied to each variable in $\mathbb X$. By Remarks \ref{remark:filter-definition} and \ref{remark:filter-consistency}, it is clear that UF is a consistent filter.

\subsection{A Consistent Bivariate Filter (BF)}

Let $(\pmb X_{1},  \dots, \pmb X_{n})$, with $\pmb X_{i} = (X_{i1}, X_{i2})^t$, be a random sample of bivariate observations. Consider also a pair of initial location and scatter estimators, 
\begin{linenomath} 
\[
\pmb T_{0n} = \left( \begin{array}{c} 
T_{0n,1} \\
T_{0n,2}
\end{array}
\right)
 \quad \text{and} \quad
\pmb C_{0n} = \left( \begin{array}{cc} 
C_{0n,11} & C_{0n,12} \\
C_{0n,21} & C_{0n,22} 
\end{array}
\right).
\]
\end{linenomath} 
Similar to the univariate case we use the coordinate-wise median and the bivariate Gnanadesikan-Kettenring estimator with MAD scale \citep{gnanadesikan:1972} for $\pmb T_{0n}$ and $\pmb C_{0n}$, respectively. More precisely, the initial scatter estimators are defined by
\begin{linenomath} 
\[
C_{0n,jk} = \frac{1}{4} \left( \text{MAD}(\{ X_{ij} + X_{ik}\})^2 - \text{MAD}(\{ X_{ij} - X_{ik}\})^2 \right), 
\]
\end{linenomath} 
where $\text{MAD}(\{Y_i\})$ denotes the MAD of $Y_1,\dots,Y_n$. Note that $C_{0n,jj} = \text{MAD}(\{X_j\})^2$, which agrees with our choice of the coordinate-wise dispersion estimators. 
Now, denote the pairwise (squared) Mahalanobis distances by $D_i = (\pmb X_{i} - \pmb T_{0n})^t \pmb C_{0n}^{-1} (\pmb X_i - \pmb T_{0n})$.
Let $G_n$ be the empirical distribution for pairwise Mahalanobis distances, 
\begin{linenomath} 
\[
G_n(t) = \frac{1}{n}\sum_{i=1}^n I( D_i \le t).
\]
\end{linenomath} 
Finally, we  filter outlying points $\pmb X_i$ by  comparing $G_n(t)$ with $G(t)$, where $G$ is a chosen reference distribution. In this paper, we use the chi-squared distribution  with two degrees of freedom, $G = \chi^2_2$.
The proportion of flagged bivariate outliers  is defined by 
\begin{linenomath} 
\begin{equation}\label{eq:2SGS-GY-d-bivariate}
\begin{aligned}
 d_{n}  &=\sup_{t\ge \eta} \left\{ G
(t)-{G}_{n}(t) \right\}^{+}. 
\end{aligned}
\end{equation}
\end{linenomath} 
Here,  $\eta = G^{-1}(\alpha)$, and we use $\alpha = 0.85$ for bivariate filtering since we now aim for moderate outliers, but other choices of $\alpha$ can be considered. Then, we flag $\lfloor nd_n \rfloor$ observations with the largest pairwise Mahalanobis distances as outlying bivariate points. Finally, the following proposition states the consistency property of the bivariate filter. 

\begin{proposition}\label{prop:2SGS-GY-asymptotic-bivariate} 
Consider a random vector $\pmb X = (X_1, X_2)^t \sim H_0$.
Also, consider a pair of bivariate  location and scatter estimators $\pmb T_{0n}$ and $\pmb C_{0n}$ such that  $\pmb T_{0n}\rightarrow\pmb\mu_0 \in \mathbb R^2$ and $\pmb C_{0n}\rightarrow\pmb\Sigma_0 \in \mathrm{PDS}(2)$ a.s. [$H_0$] $\mathrm{(}\mathrm{PDS}(q)$ is the set of all positive definite symmetric matrices of size $q\mathrm{)}$. Let $G_0(t) = P_{H_0}( (\pmb X - \pmb \mu_0)^t \pmb\Sigma_0^{-1} (\pmb X - \pmb \mu_0) \le t)$ and suppose that $G_0$ is continuous. 
If  the reference distribution $G$  satisfies:
\begin{linenomath} 
\begin{equation}\label{eq:2SGS-reference-bivariate}
\max_{t\ge \eta}\left\{G(t)-G_{0}( t) \right\}\leq 0,
\end{equation}
\end{linenomath} 
then
\begin{linenomath} 
\[
\frac{n_{0}}{n}\rightarrow0\text{ a.s.,}%
\]
\end{linenomath} 
where 
\begin{linenomath} 
\[
n_{0}= \lfloor nd_{n} \rfloor.
\]
\end{linenomath} 
\end{proposition}

In the next section, we will define the global univariate-and-bivariate filter, UBF, using UF and BF as building blocks.

\subsection{A Consistent Univariate and Bivariate Filter (UBF)}

We first apply the univariate filter from \citet{agostinelli:2014} to each variable in $\mathbb X$ separately using the initial location and dispersion estimators, $\pmb T_{0n} = (T_{0n,1}, \dots, T_{0n,p})$ and $\pmb S_{0n} = (S_{0n,1}, \dots, S_{0n,p})$.
Let $\mathbb U$ be the resulting auxiliary matrix of zeros and ones with zeros indicating the filtered entries in $\mathbb X$. We next iterate over all pairs of variables in $\mathbb X$ to identify outlying bivariate points which helps filtering the moderately contaminated cells.

 Fix a pair of variables, $ (X_{ij}, X_{ik})$ and set $\pmb X_i^{(jk)} = (X_{ij}, X_{ik})$.  
Let $\pmb C_{0n}^{(jk)}$ be an  initial pairwise scatter matrix estimator for this pair of variables, 
for example, the Gnanadesikan-Kettenring estimator. Note that pairwise scatter matrices do not ensure positive definiteness of $\pmb C_{0n}$, but this is not necessary in this case because only bivariate scatter matrix, $\pmb C_{0n}^{(jk)}$, is required 
in each bivariate filtering. 
We calculate the pairwise Mahalanobis distances $D_i^{(jk)} = (\pmb X_i^{(jk)} - \pmb T_{0n}^{(jk)})^t (\pmb C_{0n}^{(jk)})^{-1} (\pmb X_i^{(jk)} - \pmb T_{0n}^{(jk)})$ and perform the bivariate filtering on the pairwise distances with no flagged components from the univariate filtering: $\{D_i^{(jk)} : U_{ij} = 1, U_{ik}=1\}$. 
We apply this procedure to  all  pairs of variables $1\le j < k \le p$. 
Let 
\begin{linenomath} 
\[
J=\left\{  (i,j,k): D_i^{(jk)} \text{ is flagged as bivariate outlier}\right\},
\]
\end{linenomath} 
be the set of triplets which identify the pairs of cells flagged by the bivariate filter in rows $i=1,...,n$.
It remains to determine which cells $(i,j)$ in row $i$ are to be flagged as cellwise outliers.   
For each cell $(i,j)$ in the data table,  $i=1,\dots,n$ and $j=1,\dots,p,$ we count the number of
flagged pairs in the $i$-th row where cell $(i,j)$ is involved:
\begin{linenomath} 
\[
m_{ij}=\#\left\{  k: (i,j,k) \in J\right\}.
\]
\end{linenomath} 
Cells with large $m_{ij}$ are likely to correspond to univariate outliers. 
Suppose that observation $X_{ij}$ is not contaminated by cellwise contamination. Then  $m_{ij}$ approximately  follows the binomial distribution,  $Bin( \sum_{k\ne j} U_{ik}, \delta)$,   under ICM, where $\delta$ is the overall proportion of cellwise outliers that were not detected by the univariate filter. 
We flag observation $X_{ij}$ if 
\begin{linenomath} 
\begin{equation}\label{eq:UBF-condition}
	m_{ij} >  c_{ij},
\end{equation}
\end{linenomath} 
where $c_{ij}$ is the 0.99-quantile of $Bin( \sum_{k\ne j} U_{ik}, \delta)$. 
In practice we obtained good results (in both simulation and real data example) using the conservative choice $\delta = 0.10$, which is adopted in this paper. 

The filter obtained as the combination of all the univariate and the bivariate filters described above is called  UBF. The following argument shows that UBF is a consistent filter. 

By Remarks \ref{remark:filter-definition} and \ref{remark:filter-consistency}, the union of all the bivariate consistent filters (from Proposition \ref{prop:2SGS-GY-asymptotic-bivariate}) is a consistent filter. Next, applying the condition described in (\ref{eq:UBF-condition}) to the union of these bivariate consistent filters yields another consistent filter. Finally, the union of this with UF results in the consistent filter, UBF.

\subsection{The DDC Filter}

Recently,  \citet{rousseeuw:2016} proposed a new procedure  to filter and impute cellwise outliers, called  {\it DetectDeviatingCells} (DDC). DDC is a  sophisticated procedure that 
uses correlations between variables to estimate the  expected value for each cell, and then flags those with an observed value that  greatly deviates from this expected value. The DDC filter exhibited 
a very good performance when used in the first step in our two-step procedure in our simulation. However, the DDC filter  is not shown to be  consistent, as needed to ensure the overall consistency of our two-step estimation procedure. 

In view of that, we propose a new filter made by intersecting UBF and DDC (denoted here as UBF-DDC). By Remarks  \ref{remark:filter-definition} and \ref{remark:filter-consistency}, UBF-DDC
is consistent. Moreover, we will show in Section \ref{sec:GSE-MCresults} and in   \ref{sec:appendix-filter-comparison} that 
UBF-DDC is very effective, yielding the best overall performances  when used as the first step in our two-step estimation procedure.

\section{Generalized Rocke S-estimators}\label{sec:GSE-GRE}

The second step of the procedure introduces robustness against casewise outliers that went undetected in the first step. 
Data that emerged from the first step has missing values that correspond to
potentially contaminated cells. To estimate the multivariate location and
scatter matrix from that data, we  use a recently developed estimator called GSE,  briefly reviewed   below.

\subsection{Review of Generalized S-estimators}\label{sec:GSE-GSE}

Related  to $\mathbb X$ denote $\mathbb U$ the auxiliary matrix of zeros and ones, with zeros indicating the corresponding missing entries. 
Let $p_{i} = p(\pmb{U}_{i})=\sum_{j=1}^{p}U_{ij}$ be the actual dimension of the observed part of $\pmb{X}_{i}$.
Given a $p$-dimensional vector of zeros and ones $\pmb u$, a  $p$-dimensional vector $\pmb{m}$ and a  $p\times p$ matrix $\pmb{A}$,
we denote by $\pmb{m}^{(\pmb{u})}$ and $\pmb{A}^{(\pmb{u})}$  the sub-vector of  $\pmb{m}$ and the sub-matrix of $\pmb{A}$, respectively,  with columns and rows corresponding to the positive entries in $\pmb{u}$.  

Define 
\begin{linenomath} 
\[
D(\pmb{x},\pmb{m},\pmb{C})=(\pmb{x}-\pmb{m})^{t}%
\pmb{C}^{-1}(\pmb{x}-\pmb{m}) 
\]
\end{linenomath} 
the squared Mahalanobis distance and 
\begin{linenomath} 
\[
D^{\ast}(\pmb{x},\pmb{m},\pmb{C})=D(\pmb{x},\pmb{m},\pmb{C}^{\ast})
\]
\end{linenomath} 
the normalized  squared Mahalanobis distances, where $\pmb{C}^{\ast}=\pmb{C}/|\pmb{C}|^{1/p},$ so $|\pmb{C}^{\ast}|=1$, and where $|A|$ is the determinant of $A$.

Let ${\pmb{\Omega}}_{0n}$ be a $p \times p$ positive definite initial estimator. 
Given the location vector $\pmb{\mu}\in\mathbb{R}^{p}$ and a $p\times p$ positive definite matrix $\pmb{\Sigma}$, we define the generalized M-scale, $s_{GS}(\pmb{\mu},\pmb{\Sigma},{\pmb{\Omega}}_{0n}, \mathbb X, \mathbb U)$, as the solution
in $s$ to the following equation:
\begin{linenomath} 
\begin{equation}\label{eq:GSE-GSE-scale}
\sum_{i=1}^{n}c_{p(\pmb{U}_{i})}\rho\left(  \frac{D^{\ast}\left(
\pmb{X}_{i}^{(\pmb{U}_{i})},\pmb{\mu}^{(\pmb{U}_{i})},\pmb{\Sigma}^{(\pmb{U}_{i})}\right)  }{s \, c_{p(\pmb{U}_{i}%
)}\,\left\vert \pmb{\Omega}_{0n}^{(\pmb{U}_{i})}\right\vert
^{1/p(\pmb{U}_{i})}}\right)  =b\sum_{i=1}^{n}c_{p(\pmb{U}_{i})}
\end{equation}
\end{linenomath} 
where $\rho(t)$ is an even, non-decreasing in $|t|$ and bounded loss function. The tuning constants  $c_{k}$, $1\leq k\leq p$, are chosen such that
\begin{linenomath} 
\begin{equation}\label{eq:GSE-GSE-tuning-constant}
E_{\Phi}\left(  \rho\left(  \dfrac{||\pmb{X}||^{2}}{c_{k}}\right)  \right)
=b,\quad\pmb{X}\sim N_{k}(\pmb{0},\pmb{I}), 
\end{equation}
\end{linenomath} 
to ensure consistency under the multivariate normal. A common choice of $\rho$ is the  Tukey's bisquare rho function, $\rho(u) =\min(1, 1 - (1 - u)^{3})$, and $b = 0.5$, as also used in this paper.

A generalized S-estimator is then defined by
\begin{linenomath} 
\begin{equation}\label{eq:GSE-GSE}
(\pmb T_{GS}, {\pmb{C}
}_{GS}) = \arg\min_{\pmb{\mu}, \pmb{\Sigma}} s_{GS}(\pmb{\mu
}, \pmb{\Sigma}, {\pmb{\Omega}}_{0n}, \mathbb X, \mathbb U)
\end{equation}
\end{linenomath} 
subject to the constraint
\begin{linenomath} 
\begin{equation}\label{eq:GSE-GSE-constraint}
s_{GS}(\pmb{\mu}, \pmb{\Sigma},
\pmb{\Sigma}, \mathbb X, \mathbb U) = 1.
\end{equation}
\end{linenomath}

\subsection{Generalized Rocke S-estimators}\label{sec:GSE-RockeGSE}

\citet{rocke:1996} showed that if the weight function $W(x) = \rho'(x)/x$ in S-estimators is non-increasing, the efficiency of the estimators tends to one when $p \to \infty$. However, this gain in efficiency is paid for by a decrease in robustness. Not surprisingly, 
the same phenomenon  has been observed for generalized S-estimators in simulation studies.  Therefore, there is a need for new generalized S-estimators with controllable efficiency/robustness trade off. 

 \citet{rocke:1996} proposed that the $\rho$ function used to compute  S-estimators should change with the dimension to prevent loss of robustness in higher dimensions. 
The Rocke-$\rho$ function is constructed based on the fact that for large $p$ the scaled squared Mahalanobis distances for normal data 
\begin{linenomath} 
\[
\frac{D(\pmb X, \pmb \mu, \pmb\Sigma)}{\sigma} \approx \frac{Z}{p} \quad \text{with} \quad Z \sim \chi^2_p,
\]
\end{linenomath} 
and hence that $D/\sigma$ are increasingly concentrated around one. So, to have a high enough, but not too high, efficiency, we should give a high weight to the values of $D/\sigma$ near one and downweight the cases where $D/\sigma$ is far from one. 

Let 
\begin{linenomath} 
\begin{equation}\label{eq:GSE-Rocke-rho-gamma}
\gamma = \min\left(\frac{\chi^2(1 - \alpha)}{p} - 1, 1 \right),
\end{equation}
\end{linenomath} 
where $\chi^2(\beta)$ is the $\beta$-quantile of $\chi^2_p$. 
In this paper, we use a conventional choice of $\alpha = 0.05$ that gives an acceptable  efficiency of the estimator. We have also explored smaller values of $\alpha$ according to \citet{maronna:2015b}, but we have seen some degree of trade-offs between efficiency and casewise robustness (see the supplementary material).
\citet{maronna:2006} proposed a modification of  the Rocke-$\rho$ function, namely 
\begin{linenomath} 
\begin{equation}\label{eq:GSE-Rocke-rho}
\rho(u) = \begin{cases}
0 & \text{for}\quad  0 \le u \le 1 - \gamma \\
\left(\frac{u-1}{4\gamma}\right)\left[3 - \left(\frac{u-1}{\gamma}\right)^2\right] + \frac{1}{2} & \text{for}\quad 1 - \gamma < u <1 + \gamma\\
1 & \text{for}\quad u \ge 1 + \gamma
\end{cases}
\end{equation}
\end{linenomath} 
which has as derivative the desired weight function that vanishes for $u \not\in [1 - \gamma, 1 + \gamma]$
\begin{linenomath} 
\[
W(u) = \frac{3}{4\gamma}\left[1 - \left(\frac{u - 1}{\gamma} \right)^2 \right] I( 1 - \gamma \le u \le 1 + \gamma).
\]
\end{linenomath} 

\begin{figure}[t!]
\centering
\includegraphics[scale=0.65]{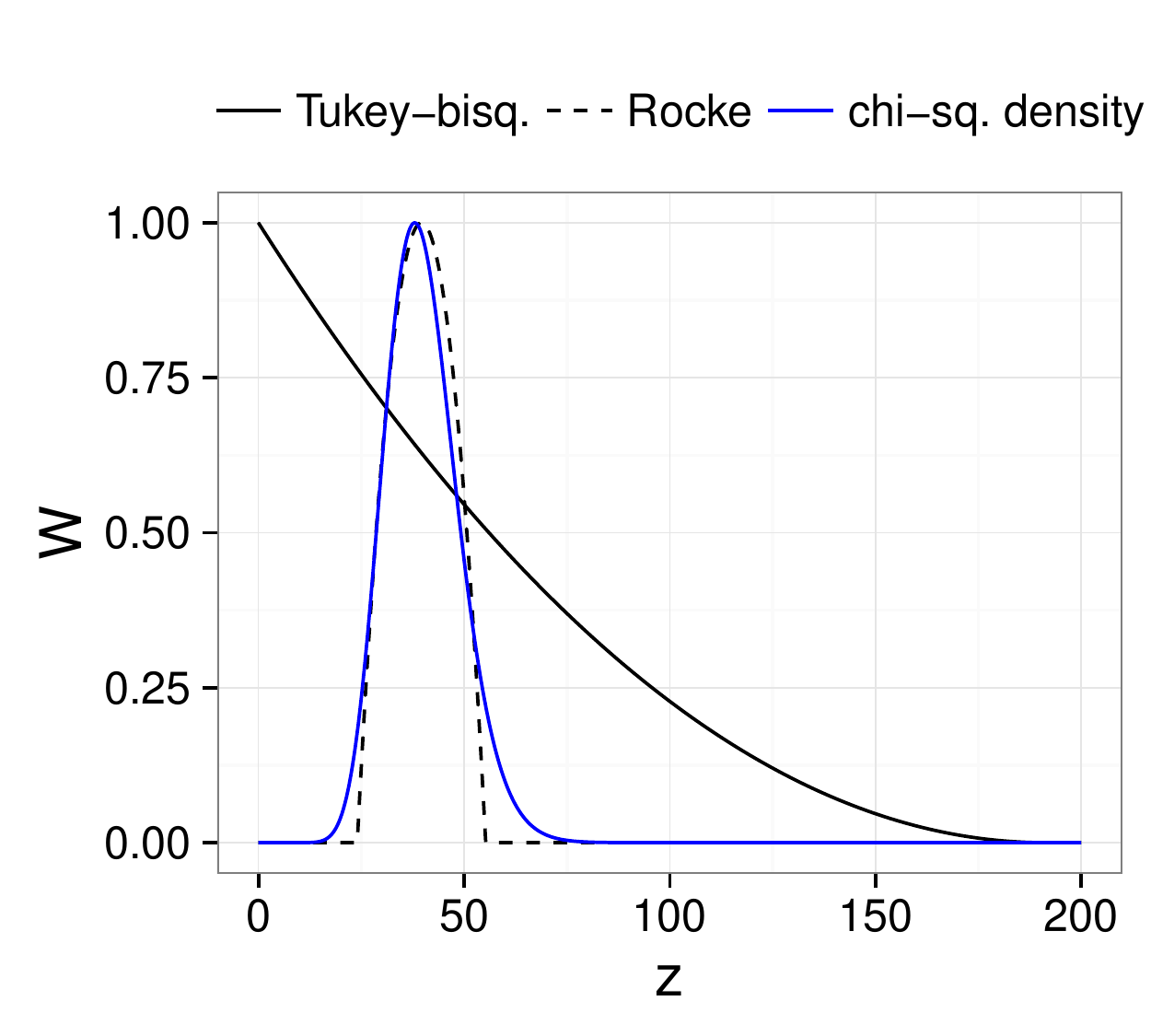} 
\caption{Weight functions of the Tukey-bisquare and the Rocke for $p=40$. Chi-square density functions are also plotted in blue for comparison. All the functions are scaled so that their maximum is 1 to facilitate comparison.}\label{fig:GSE-example-weights}
\end{figure}

Figure \ref{fig:GSE-example-weights} compares the Rocke-weight function, $W_{Rocke}(z/c_p)$, and the Tukey-bisquare weight function, $W_{Tukey}(z/c_p)$, for $p=40$, where $c_p$ as defined in (\ref{eq:GSE-GSE-tuning-constant}). The chi-square density function is also plotted in blue for comparison. When $p$ is large the tail of the Tukey-bisquare weight function greatly deviates from the tail of the chi-square density function and inappropriately assigns high weights to large distances. On the other hand, the Rocke-weight function can resemble the shape of the chi-square density function and is capable of assigning low weights to large distances. 

Finally, we define the generalized Rocke S-estimators or GRE by (\ref{eq:GSE-GSE}) and (\ref{eq:GSE-GSE-constraint}) with  the $\rho$-function in (\ref{eq:GSE-GSE-scale}) replaced by the modified Rocke-$\rho$ function in (\ref{eq:GSE-Rocke-rho}). We compared GRE with GSE via  simulation  and found that  GRE has a substantial better performance in dealing with casewise outliers when $p$ is large (e.g., $p > 10$). Results from this simulation study are provided in the  supplementary material.

\section{Computational  Issues}\label{sec:GSE-computing-issue}

The generalized S-estimators described above are computed via iterative re-weighted means and covariances, starting from an initial estimate. We now discuss some  computing issues associated with this iterative  procedure.

\subsection{Computation of the Initial Estimator}

For the initial estimate, the extended minimum volume ellipsoid (EMVE) has been used, as suggested by \citet[][]{danilov:2012}. The EMVE is computed with a large number of subsamples ($>500$) to increase the chance that at least one clean subsample is obtained. 
 Let $\varepsilon$ be the proportion of contamination in the data and $m$ be the subsample size. The probability of having at least one clean subsample of size $m$ out of $M$ subsamples is
\begin{linenomath} 
\begin{equation}\label{eq:GSE-nresample-calculation}
q = 1 - \left[1 - \left(\begin{array}{c} n \cdot(1 - \varepsilon) \\ m \end{array}\right)/\left(\begin{array}{c} n \\ m \end{array}\right) \right]^M.
\end{equation}
\end{linenomath} 
For large $p$, the number of subsamples $M$ required for a large $q$, say $q = 0.99$, can be impractically large, dramatically slowing down the computation. For example, suppose $m = p$, $n = 10p$, and $\varepsilon = 0.50$. If $p = 10$, then $M = 7758$; if $p = 30$, then $M = 2.48 \times 10^{10}$; and if $p = 50$, then $M = 4.15 \times 10^{16}$. Therefore, there is a need for a faster and more reliable starting point for large $p$.

Alternatively,  pairwise scatter estimators could be used as fast initial estimator \citep[e.g.,][]{alqallaf:2002}. Previous simulation studies have shown that pairwise scatter estimators are robust against cellwise outliers, but they perform not as well in the presence of casewise outliers and finely shaped multivariate data \citep{danilov:2012, agostinelli:2014}.

\subsubsection{Cluster-Based Subsampling}
Next, we introduce a cluster-based algorithm for faster and more reliable subsampling for the computation of EMVE. The EMVE computed with the cluster-based subsampling is called called EMVE-C  throughout the paper.

High-dimensional data have several interesting geometrical properties as  described  in \citet{hall:2005}. One such property that motivated the Rocke-$\rho$ function, as well as the following algorithm, is that for large $p$ the $p$-variate standard normal distribution $N_p(\pmb 0, \pmb I)$ is concentrated ``near" the spherical shell with radius $\sqrt p$. So, if outliers have a slightly different covariance structure from  clean data, they would appear geometrically different. Therefore, we could apply a clustering algorithm to first separate the outliers from the clean data. Subsampling from a big cluster, which in principle is composed  of mostly clean cases, should be more reliable and require fewer number of subsamples.

Given $\mathbb X$ and $\mathbb U$. The following steps describe our clustering-based subsampling:  
\begin{enumerate}\itemsep0pt
\item Standardize the data $\mathbb X$ with some initial location and dispersion estimator $T_{0j}$ and $S_{0j}$. Common choices for $T_{0j}$ and $S_{0j}$ that are also adopted in this paper are the coordinate-wise median and MAD. Denote the standardized data by $\mathbb Z = (\pmb Z_1, \dots, \pmb Z_n)^t$, where $\pmb Z_i = (Z_{i1}, \dots, Z_{ip})^t$ and $Z_{ij} = (X_{ij} - T_{0j})/S_{0j}$.

\item Compute a simple robust correlation matrix estimate $\pmb R = (R_{jk})$. Here, we use the  Gnanadesikan-Kettenring estimator \citep{gnanadesikan:1972}, where
\begin{linenomath} 
\[
R_{ij} = \frac{1}{4} (S_{0jk+}^2 - S_{0jk-}^2),
\]
\end{linenomath} 
and where $S_{0jk+}$ is the dispersion estimate for $\{ Z_{ij} + Z_{ik} | U_{ij} = 1, U_{ik} = 1 \}$ and $S_{0jk-}$ the estimate for $\{ Z_{ij} - Z_{ik} | U_{ij} = 1, U_{ik} = 1\}$.  We use $Q_n$ \citep{rousseeuw:1993} for the dispersion estimate. 

\item Compute the eigenvalues $\lambda_1 \ge \dots \ge \lambda_p$ and eigenvectors $\pmb e_1, \dots, \pmb e_p$ of the correlation matrix estimate 
\begin{linenomath} 
\[
\pmb R = \pmb E \pmb \Lambda \pmb E^t, 
\]
\end{linenomath} 
where $\pmb\Lambda = \text{diag}(\lambda_1, \dots, \lambda_p)$ and $\pmb E = (\pmb e_1, \dots, \pmb e_p)$. Let $p_+$ be the largest dimension such that $\lambda_{j} > 0$ for $j =1, \dots, p_+$. Retain only the eigenvectors $\pmb E_0 = (\pmb e_1, \dots, \pmb e_{p_+})$ with a positive eigenvalue.

\item Complete the standardized data $\mathbb Z$ by replacing each missing entry, as indicated by $\mathbb U$, by zero. Then, project the data onto the basis eigenvectors $\tilde{\pmb Z} = \pmb Z \pmb E_0$, and then standardize the columns of $\tilde{\pmb Z}$, or so called principal components, using coordinate-wise median and MAD of $\tilde{\pmb Z}$.

\item Search for a ``clean" cluster $C$ in the standardized $\tilde{\pmb Z}$ using a hierarchical clustering framework by doing the following. First, compute the dissimilarity matrix for the principal components using the  Euclidean metric. Then, apply classical hierarchical clustering (with any linkage of choice). A common choice is the Ward's linkage, which is  adopted in this paper. Finally, define the ``clean" cluster by the smallest sub-cluster $C$ with a size at least $n/2$. This can be obtained by cutting the clustering tree at various heights from the top until all the clusters have size less than $n/2$. 

\item Take a subsample of size $n_0$ from $C$. 
\end{enumerate}

With good clustering results, we can draw fewer subsamples, and equally important, we can use a larger subsample size. The current default  choices in GSE  are $M=500$ subsamples of size $n_0 = (p + 1)/(1 - \alpha_{mis})$ as suggested in \citet{danilov:2012}, where $\alpha_{mis}$ is the fraction of missing data ($\alpha_{mis} =$ number of missing entries $/(np)$). For the new clustering-based subsampling, we choose $M=50$ and $n_0 = 2(p + 1)/(1 - \alpha_{mis})$  in view of  their overall good  performance in our simulation study. However, using equation (\ref{eq:GSE-nresample-calculation}), a more formal procedure for the choice of  $M$ and $n_0$ could be considered. 
 $M$ and $n_0$ could be chosen as a function of the cluster size $C$, the expected remaining fraction of contamination $\delta$, and a desired level of confidence. In such case,    $n$ and $\varepsilon$ in equation (\ref{eq:GSE-nresample-calculation})  should  be replaced by  to the size of the cluster $C$ and the value of $\delta$, respectively. 
Without clustering, $\varepsilon$ would be  chosen fairly large (e.g.  $\varepsilon =0.50$) for conservative reasons. However,  with clustering,  $\varepsilon$ can be made smaller (e.g., $\varepsilon \le  0.10$). 

In general, $p$ is the primary driver of computational time, but the procedure could also be time-consuming for large $n$ because the number of operations required by hierarchical clustering is of order $n^3$. As an alternative, one may bypass the hierarchical clustering step and sample directly from the data points with the smallest Euclidean distances to the origin calculated from $\tilde{\pmb Z}$. This is because the Euclidean distances, in principle, should approximate the Mahalanobis distances to the mean of the original data.  However, our simulations show that the hierarchical clustering step is essential for the excellent performance of the estimates, and that this step entails only a small increase in real computational time, even for large $n$. 

A recent simulation study \citep{maronna:2015b} has shown that Rocke estimator starting from the  the ``kurtosis plus specific direction" (KSD) estimator \citep{pena:2001}  estimator  can attain high efficiency and high robustness for large  $p$. The KSD estimator uses a multivariate outlier detection procedure based on finding directions that maximize or minimize the kurtosis coefficient of the respective projections. The ``clean" cases that were not flagged as outliers are then used for estimating multivariate location and scatter matrix. Unfortunately, KSD is not implemented for incomplete data. The study of the adaption of KSD for incomplete data would be of interest and worth of future research.

\subsection{Other  Computational  Issues}

There is no formal proof that the recursive algorithm  decreases the objective function at each iteration for the case of  generalized S-estimators with a monotonic weight function  \citep{danilov:2012}. This also the case for generalized S-estimators with a non-monotonic weight function. For Rocke estimators with complete data,  \citet[][see Section 9.6.3]{maronna:2006} described an algorithm that ensures attaining a local minimum. We have adapted this algorithm for the generalized counterparts. Although we cannot provide a formal proof,  we have seen so far in our experiments that the descending property of the recursive algorithms always holds.

\section{Two-Step Estimation and Simulation Results}\label{sec:GSE-MCresults}

The original two-step approach for global--robust estimation under cellwise and casewise contamination  is to first flag outlying cells in the data table and to replace them by NA's using a univariate filter only (shortened to UF). In the second step, the generalized S-estimator is then applied to this incomplete data. Our new version of this is to replace UF in the first step by the proposed combination of univariate-and-bivariate filter and DDC (shortened to UBF-DDC) and to replace GSE in the second step by GRE-C (i.e., GRE starting from EMVE-C).  We call the new two-step procedure UBF-DDC-GRE-C. The new procedure will be made available in the \texttt{TSGS} function
in the \texttt{R} package \texttt{GSE} \citep{leung:2015}.

We now conduct a simulation study similar to that in \citet{agostinelli:2014} to compare  the two-step procedures, UF-GSE as introduced in \citet{agostinelli:2014} and UBF-DDC-GRE-C, as well as 
the classical correlation estimator (MLE) and several other robust estimators that showed  a competitive performance under 
\begin{itemize}
\item Cellwise contamination: SnipEM (shortened to Snip) introduced in \citet{farcomeni:2014a} 
\item Casewise contamination:  Rocke S-estimator as recently revisited by \citet{maronna:2015b} and HSD introduced by \citet{vanAelst:2012}
\item Cellwise  and casewise contamination: DetMCDScore (shortened to DMCDSc)  introduced by \citet{rousseeuw:2015}
\end{itemize}
We also considered the different variations of the two-step procedures using different first steps, including UBF-GRE-C and DDC-GRE-C. However, UBF-DDC-GRE-C generally performs better in simulations than UBF-GRE-C and DDC-GRE-C. Therefore, we present only the results of UBF-DDC-GRE-C here. The complete results of UBF-GRE-C and DDC-GRE-C can be found in  \ref{sec:appendix-filter-comparison}.

We consider clean and contaminated samples from a $N_p(\pmb{\mu_0},\pmb{\Sigma_0})$ distribution with dimension $p=10,20,30,40,50$  and sample size $n= 10p$. The simulation mechanisms are briefly described below. 

Since the contamination models and  the estimators  considered in our simulation study
 are location and scale equivariant, we can assume without loss of
generality that the mean, $\pmb\mu_0$, is equal to $\pmb 0$  and the variances in  $\text{diag}(\pmb\Sigma_0)$ are all equal to $\pmb 1$. That is, $\pmb\Sigma_0$ is a correlation matrix. 

Since the cellwise contamination model and the estimators are not affine-equivariant, we consider the two different approaches to introduce correlation structures: 
\begin{itemize}
\item Random correlation as described in \citet{agostinelli:2014} and 
\item First order autoregressive correlation.
\end{itemize}
The random correlation structure generally has small correlations, especially with increasing $p$. For example, for $p=10$, the maximum correlation values have an average of $0.49$, and for $p=50$, the average maximum is $0.28$. 
So, we consider the first order autoregressive correlation (AR1) with higher correlations, in which the correlation matrix has entries
\begin{linenomath}
\[
\Sigma_{0,jk} = \rho^{|j-k|},
\]
\end{linenomath}
with $\rho=0.9$. 

We then consider the following scenarios:
\begin{itemize}
\item Clean data: No further changes are done to the data.
\item Cellwise contamination: We randomly replace a $\epsilon$  of the cells in the data matrix by $X_{ij}^{cont} \sim N(k, 0.1^2)$, where $k=1,2,\dots,10$.
\item Casewise contamination: We randomly replace a $\epsilon$ of the cases in the data matrix by $\pmb X_i^{cont} \sim 0.5  N( c  \pmb v,  0.1^2  \pmb I) + 0.5  N(- c  \pmb v,  0.1^2  \pmb I)$, where $c=\sqrt{k (\chi^2)^{-1}_p(0.99)}$ and $k = 1,2,\dots, 20$ and $\pmb{v}$ is the eigenvector corresponding to the smallest eigenvalue of $\pmb{\Sigma}_{0}$ with length such that $\left(\pmb{v} -\pmb\mu_0\right)^{t}\pmb{\Sigma}_{0}^{-1}
\left(\pmb{v}-\pmb\mu_0\right)=1$. Experiments show that the placement of outliers in this way is the least favorable for the proposed estimator. 
\end{itemize}
We consider $\epsilon =0.02, 0.05$ for cellwise contamination, and $\epsilon = 0.10, 0.20$ for casewise contamination. 
The number of replicates in our simulation study is $N=500$. 

The performance of a given scatter estimator $\pmb{\Sigma}_n$ is
measured by the Kulback--Leibler  divergence between two Gaussian distribution with the same mean and 
covariances $\pmb{\Sigma}$ and $\pmb{\Sigma}_{0}$:
\begin{linenomath}
\[
D(\pmb{\Sigma}, \pmb{\Sigma}_{0}) = \mbox{trace}(\pmb{\Sigma
}\pmb{\Sigma}_{0}^{-1}) - \log(|\pmb{\Sigma}\pmb{\Sigma}_{0}^{-1}|) -
p.
\]
\end{linenomath}
This divergence also appears  in the likelihood ratio test statistics for testing 
the null hypothesis that a multivariate normal distribution has covariance
matrix $\pmb{\Sigma}= \pmb{\Sigma}_{0}$. We call this divergence measure the likelihood ratio test distance (LRT). Then, the performance of an
estimator $\pmb{\Sigma}_n$ is summarized by
\begin{linenomath}
\[
\overline{D}({\pmb{\Sigma}}_n, \pmb{\Sigma}_{0}) =\frac{1}{N} \sum_{i=1}^{N}
D(\hat{\pmb{\Sigma}}_{n,i}, \pmb{\Sigma}_{0})
\]
\end{linenomath}
where $\hat{\pmb{\Sigma}}_{n,i}$ is the estimate at the $i$-th replication. Finally, the maximum average LRT distances over all considered contamination values, $k$, is also calculated.

\begin{table}[t!]
\centering
\footnotesize
\caption{Maximum average LRT distances under cellwise contamination. The sample size is $n=10p$. }\label{tab:GSE-2SGRE-Cellwise}
\begin{tabular}{llcccccccHHHHc}
  \hline
Corr. & $p$ & $\epsilon$ & MLE & Rocke & HSD & Snip & DMCDSc & UF- & UBF- & UF- & UBF- & DDC- & UBF-DDC- \\ 
& & & & & & &  & GSE & GSE & GRE-C & GRE-C  & GRE-C & GRE-C\\
\hline
 Random & 10 & 0 & 0.6 & 1.2 & 0.8 & 5.0 & 1.5 & 0.8 & 0.9 & 1.2 & 1.3 &  1.0 & 1.0\\
   &  & 0.02 & 114.8 & 1.2 & 2.3 & 6.9 & 1.6 & 1.2 & 1.4 & 1.3 & 1.4 & 1.1 & 1.1\\ 
   &  & 0.05 & 285.4 & 3.6 & 11.2 & 7.5 & 3.2 & 4.5 & 4.4 & 2.2 & 2.5 & 2.6 & 2.5\\ 
   & 20 & 0 & 1.1 & 2.0 & 1.2 & 11.5 & 2.0 & 1.3 & 1.5 & 1.9 & 2.0 & 1.8 & 1.8 \\ 
   &  & 0.02 & 146.1 & 2.7 & 10.6 & 13.9 & 2.6 & 4.0 & 4.4 & 2.9 & 3.0 & 2.5 & 2.5  \\ 
   &  & 0.05 & 375.9 & 187.2 & 57.1 & 15.5 & 9.3 & 11.0 & 11.1 & 8.0 & 8.2 & 7.7 & 7.3\\ 
   & 30 & 0 & 1.6 & 2.8 & 1.7 & 16.7 & 2.6 & 1.9 & 2.0 & 3.4 & 3.9 & 3.5 & 3.3\\ 
   &  & 0.02 & 179.0 & 23.1 & 22.6 & 18.5 & 4.4 & 5.8 & 6.3 & 5.4 & 5.9 & 5.3 & 5.0 \\ 
   &  & 0.05 & 475 & 380.5 & 123.1 & 20.8 & 13.7 & 14.2 & 14.8 & 12.3 & 13.4 & 14.2 & 13.3 \\  
   & 40 & 0 & 2.1 & 3.6 & 2.3 & 20.7 & 3.2 & 2.4 & 2.6 & 5.9 & 6.2 & 5.8 & 5.8 \\ 
   &  & 0.02 & 215.1 & 121.3 & 38.9 & 22.6 & 6.0 & 7.3 & 8.0 & 9.4 & 10.9 & 9.5 & 8.8 \\ 
   &  & 0.05 & $>$500 & $>$500 & 212.4 & 25.8 & 17.9 & 16.6 & 17.4 & 18.4 & 19.9 & 18.8 & 18.6 \\ 
   & 50 & 0 & 2.7 & 4.4 & 2.8 & 25.4 & 3.8 & 2.9 & 3.2 & 5.2 & 5.3  & 4.9 & 4.9 \\ 
   &  & 0.02 & 249.0 & 192.8 & 58.7 & 27.1 & 8.1 & 9.1 & 10.0 & 12.5 & 12.9 & 12.5 & 12.1 \\ 
   &  & 0.05 & $>$500 & $>$500 & 298.7 & 29.7 & 20.7 & 19.6 & 20.6 & 22.7 & 23.6 & 24.4 & 23.8 \\ 
   \hline
AR1(0.9) & 10 & 0 & 0.6 & 1.1 & 0.8 & 4.3 & 1.4 & 0.7 & 0.8 & 1.1 & 1.2 & 1.1 & 1.0 \\ 
   &  & 0.02 & 149.8 & 1.2 & 0.9 & 4.9 & 1.5 & 0.9 & 0.9 & 1.2 & 1.3 & 1.1 & 1.0 \\ 
   &  & 0.05 & 383.8 & 2.6 & 2.8 & 7.0 & 3.1 & 2.1 & 1.1 & 1.7 & 1.4 & 1.3 & 1.3 \\ 
   & 20 & 0 & 1.1 & 1.9 & 1.2 & 7.8 & 2.1 & 1.2 & 1.3 & 1.8 & 1.9 & 1.8 & 1.7 \\ 
   &  & 0.02 & 311.3 & 2.5 & 3.9 & 10.5 & 2.6 & 2.1 & 1.5 & 2.2 & 2.1 & 2.0 & 1.9 \\ 
   &  & 0.05 & $>$500 & $>$500 & 31.3 & 14.3 & 12.3 & 9.3 & 2.7 & 7.6 & 2.8 & 2.1 & 2.5 \\ 
   & 30 & 0 & 1.6 & 2.8 & 1.8 & 9.4 & 2.7 & 1.7 & 1.8 & 3.2 & 3.4 & 3.6 & 3.2 \\ 
   &  & 0.02 & 475.9 & 71.1 & 10.7 & 13.9 & 5.4 & 4.0 & 2.3 & 3.9 & 3.4 & 3.5 & 3.3 \\ 
   &  & 0.05 & $>$500 & $>$500 & 103.3 & 19.8 & 22.6 & 20.3 & 6.2 & 18.1 & 5.5 & 3.4 & 3.6 \\ 
   & 40 & 0 & 2.1 & 3.6 & 2.2 & 10.9 & 3.4 & 2.3 & 2.3 & 5.5 & 5.7 &  5.8 & 5.5 \\ 
   &  & 0.02 & $>$500 & 222.1 & 22.7 & 16.2 & 8.9 & 6.7 & 3.5 & 6.5 & 5.7 & 6.0 & 5.6 \\ 
   &  & 0.05 & $>$500 & $>$500 & 259.9 & 23.7 & 34.8 & 31.4 & 14.0 & 29.7 & 12.4 & 6.1 & 5.9 \\ 
   & 50 & 0 & 2.7 & 4.4 & 2.8 & 13.0 & 4.0 & 2.8 & 2.9 & 5.5 & 5.2 & 4.6 & 5.0 \\  
   &  & 0.02 & $>$500 & $>$500 & 43.3 & 18.9 & 12.8 & 9.7 & 4.9 & 9.7 & 6.4 & 6.4 & 7.8 \\ 
   &  & 0.05 & $>$500 & $>$500 & $>$500 & 28.9 & 46.5 & 42.8 & 22.6 & 40.8 & 20.4 & 7.9 & 8.9 \\ 
   \hline
\end{tabular}
\end{table}

\begin{figure}[t!]
\centering
\includegraphics[scale=0.5]{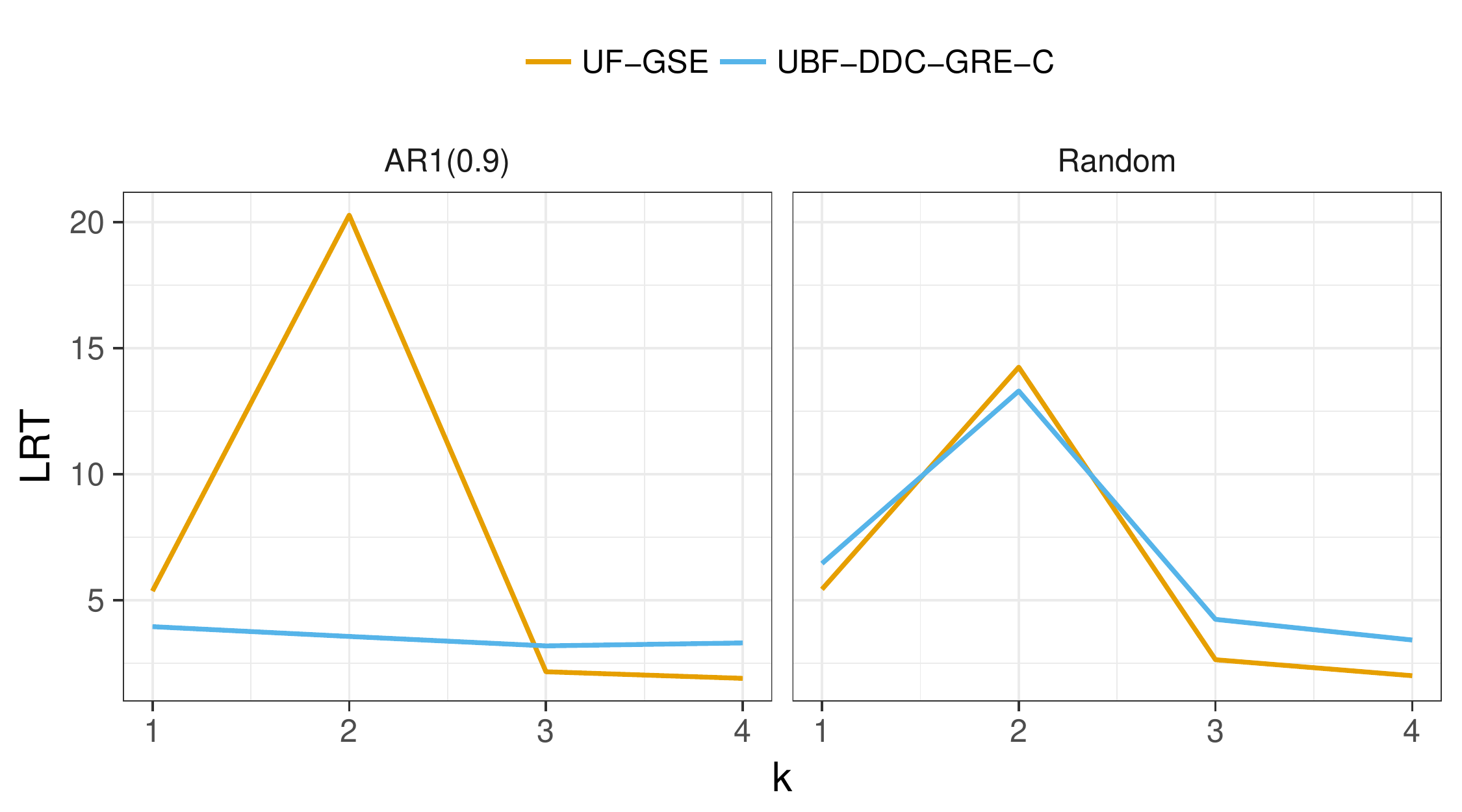}
\caption{Average LRT distance behaviors for various contamination values, $k$, of UF-GSE and  UBF-DDC-GSE for random and AR$1(0.9)$ correlations under 5\% cellwise contamination. The dimension is $p=30$ and the sample size is $n=10p$. The results remain the same for larger values of $k$; thus, they are not included in the figure. }\label{fig:2SGR-LRT-cellwise-curve}
\end{figure}

Table \ref{tab:GSE-2SGRE-Cellwise} shows the maximum average LRT distances  under cellwise contamination. 
UBF-DDC-GRE-C and UF-GSE perform similarly under random correlation, but UBF-DDC-GRE-C outperforms UF-GSE under AR$1(0.9)$. When correlations are small, like in  random correlation,  the bivariate filter fails to filter moderate cellwise outliers (e.g., $k=2$) because there is not enough information about the bivariate correlation structure in the data. Therefore, the bivariate filter gives similar results as the univariate filter. However, when correlations are large, like in AR$1(0.9)$, the bivariate filter can filter moderate cellwise outliers and therefore, outperforms the univariate filter. This is  demonstrated, for example, in Figure \ref{fig:2SGR-LRT-cellwise-curve}  which shows the average LRT distance behaviors for various cellwise contamination values, $k$.

\begin{table}[t!]
\centering
\footnotesize
\caption{Maximum average LRT distances under casewise contamination. The sample size is $n=10p$. }\label{tab:GSE-2SGRE-Casewise}
\begin{tabular}{llcccccccHHHHc}
  \hline
Corr. & $p$ & $\epsilon$ & MLE & Rocke & HSD & Snip & DMCDSc & UF- & UBF- & UF- & UBF- & DDC- & UBF-DDC- \\ 
& & & & & & &  & GSE & GSE & GRE-C & GRE-C  & GRE-C & GRE-C\\
   \hline
 Random & 10 & 0 & 0.6 & 1.2 & 0.8 & 5.0 & 1.5 & 0.8 & 0.9 & 1.2 & 1.3 & 1.0 & 1.0 \\
   &  & 0.10 & 43.1 & 2.8 & 3.9 & 44.4 & 4.9 & 9.7 & 18.5 & 11.0 & 19.1 & 9.4 & 7.7  \\ 
   &  & 0.20 & 89.0 & 4.7 & 21.8 & 110.3 & 123.6 & 91.8 & 146.8 & 30.1 & 53.0 & 25.3 & 23.7 \\ 
   & 20 & 0 & 1.1 & 2.0 & 1.2 & 11.5 & 2.0 & 1.3 & 1.5 & 1.9 & 2.0 & 1.8 & 1.8 \\ 
   &  & 0.10 & 77.0 & 3.4 & 13.4 & 76.9 & 37.8 & 29.7 & 50.1 & 11.5 & 20.9 & 9.5 & 9.1  \\ 
   &  & 0.20 & 146.7 & 5.6 & 95.9 & 166.5 & 187.6 & 291.8 & 311.4 & 22.0 & 49.3 &  18.0 & 17.4 \\ 
   & 30 & 0 & 1.6 & 2.8 & 1.7 & 16.7 & 2.6 & 1.9 & 2.0 & 3.4 & 3.9 & 3.5 & 3.3 \\ 
   &  & 0.10 & 100.0 & 4.3 & 26.1 & 82.3 & 118.6 & 75.3 & 101.3 & 12.8 & 21.8 & 10.6 & 9.9 \\  
   &  & 0.20 & 200.7 & 7.4 & 297.7 & 220.9 & 268.4 & 415.5 & 445.2 & 21.7 & 47.6 &  18.7 & 16.9 \\ 
   & 40 & 0 & 2.1 & 3.6 & 2.3 & 20.7 & 3.2 & 2.4 & 2.6 & 5.9 & 6.2 & 5.8 & 5.8 \\ 
   &  & 0.10 & 125.9 & 5.2 & 46.3 & 101.6 & 130.6 & 140.2 & 168.8 & 18.6 & 29.5 & 17.7 & 16.2 \\ 
   &  & 0.20 & 252.4 & 9.1 & $>$500 & 186.2 & 340.1 & $>$500 & 579.9 & 22.7 & 52.3 & 21.2 & 19.5 \\ 
   & 50 & 0 & 2.7 & 4.4 & 2.8 & 25.4 & 3.8 & 2.9 & 3.2 & 5.2 & 5.3 & 4.9 & 4.9 \\ 
   &  & 0.10 & 150.3 & 5.9 & 80.0 & 121.9 & 139.5 & 258.1 & 228.8 & 27.5 & 43.4 & 21.2 & 17.6 \\ 
   &  & 0.20 & 303.1 & 10.0 & $>$500 & 224.3 & 407.7 & $>$500 & $>$500 & 24.2 & 64.8 & 23.7 & 23.0 \\ 
   \hline
AR1(0.9) & 10 & 0 & 0.6 & 1.1 & 0.8 & 4.3 & 1.4 & 0.7 & 0.8 & 1.1 & 1.2 & 1.1 & 1.0 \\ 
   &  & 0.10 & 43.1 & 2.8 & 1.7 & 20.2 & 2.9 & 3.7 & 4.3 & 3.1 & 3.6 & 3.0 & 2.9 \\ 
   &  & 0.20 & 88.9 & 4.8 & 8.7 & 49.7 & 29.7 & 50.8 & 50.1 & 7.2 & 8.4 & 6.8 & 6.9 \\ 
   & 20 & 0 & 1.1 & 1.9 & 1.2 & 7.8 & 2.1 & 1.2 & 1.3 & 1.8 & 1.9 & 1.8 & 1.7 \\ 
   &  & 0.10 & 77.0 & 2.8 & 4.7 & 43.8 & 14.8 & 12.9 & 14.9 & 3.5 & 4.3 & 3.3 & 3.3 \\ 
   &  & 0.20 & 146.6 & 5.3 & 35.3 & 113.0 & 87.6 & 260.5 & 193.9 & 7.3 & 10.5 & 6.0 & 6.0\\ 
   & 30 & 0 & 1.6 & 2.8 & 1.8 & 9.4 & 2.7 & 1.7 & 1.8 & 3.2 & 3.4 & 3.6 & 3.2 \\ 
   &  & 0.10 & 98.9 & 3.4 & 8.9 & 66.1 & 32.2 & 31.3 & 37.7 & 4.1 & 5.1 & 4.2 & 4.1 \\ 
   &  & 0.20 & 200.5 & 8.2 & 155.5 & 144.8 & 122.9 & 372.7 & 365.1 & 8.4 & 13.3 & 6.9 & 6.8 \\ 
   & 40 & 0 & 2.1 & 3.6 & 2.2 & 10.9 & 3.4 & 2.3 & 2.3 & 5.5 & 5.7  & 5.8 & 5.5 \\ 
   &  & 0.10 & 124.9 & 4.3 & 15.6 & 83.7 & 49.2 & 69.1 & 75.5 & 6.4 & 7.3 & 5.8 & 6.4 \\ 
   &  & 0.20 & 253.0 & 9.2 & 430.3 & 151.9 & 209.3 & 477.6 & 479.7 & 10.0 & 17.4 & 8.9 & 8.7 \\ 
   & 50 & 0 & 2.7 & 4.4 & 2.8 & 13.0 & 4.0 & 2.8 & 2.9 & 5.5 & 5.2 & 4.6 & 5.0 \\ 
   &  & 0.10 & 150.2 & 5.1 & 26.5 & 103.3 & 64.4 & 148.2 & 160.1 & 7.6 & 8.1 & 7.5 & 7.9 \\ 
   &  & 0.20 & 302.6 & 10.1 & $>$500 & 188.5 & 276.0 & $>$500 & $>$500 & 11.0 & 21.2 & 10.0 & 8.8 \\ 
   \hline
\end{tabular}
\end{table}

\begin{figure}[t!]
\centering
\includegraphics[scale=0.5]{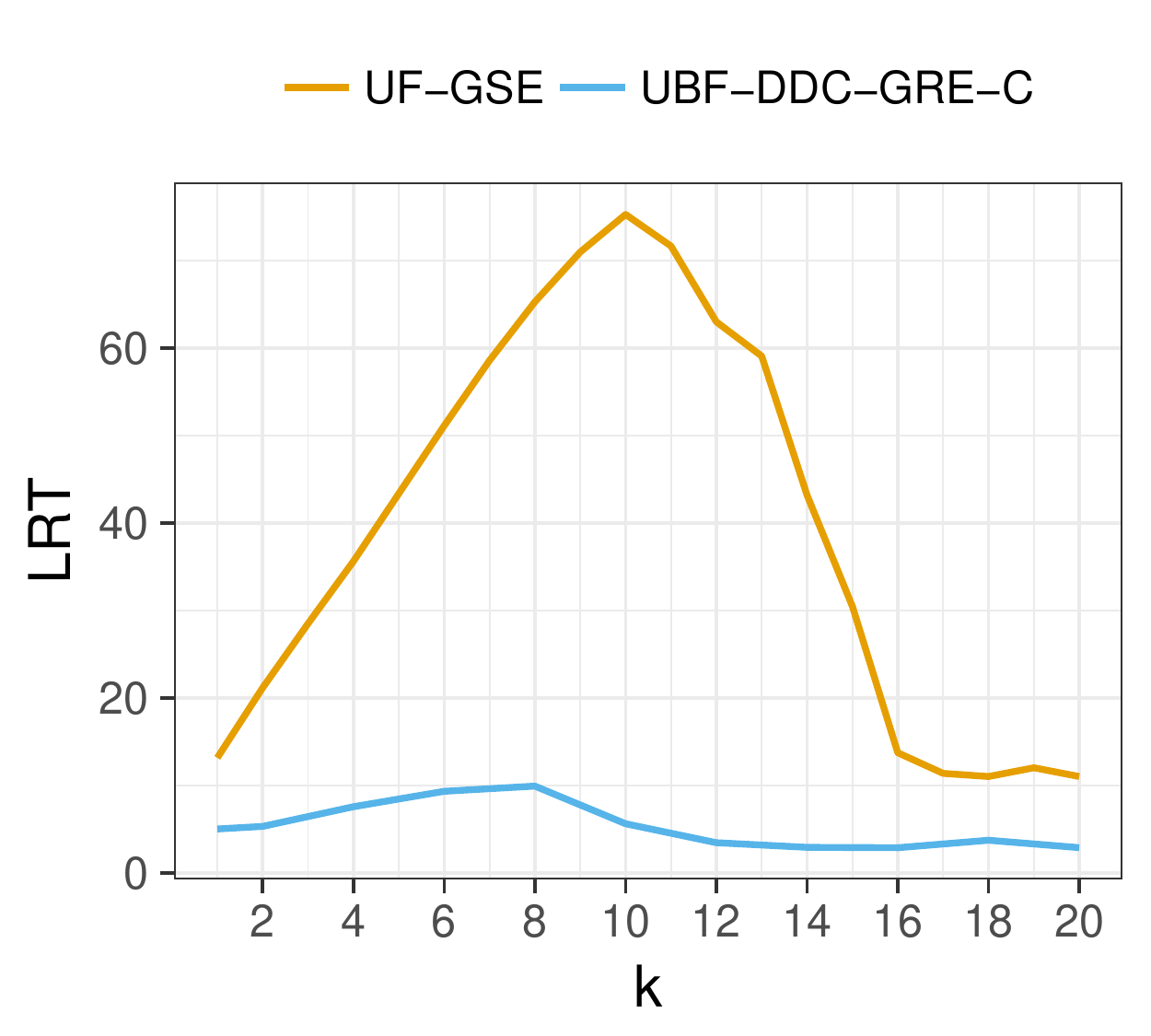}
\caption{Average LRT distance behaviors for various contamination values, $k$, of UF-GSE and UBF-DDC-GRE-C for random correlations under 10\% casewise contamination. The dimension is $p=30$ and the sample size is $n=10p$.  }\label{fig:GSE-2SGR-LRT-casewise-curve}
\end{figure}

Table \ref{tab:GSE-2SGRE-Casewise} shows the maximum average LRT distances  under casewise contamination. Overall, UBF-DDC-GRE-C outperforms UF-GSE. This is because the Rocke $\rho$ function in GRE in UBF-DDC-GRE-C is more capable of downweighting moderate casewise outliers (e.g., $10 < k < 20$) than the Tukey-bisquare $\rho$ function in GSE in UF-GSE. Therefore, UBF-DDC-GRE-C outperforms UF-GSE under moderate casewise contamination and gives overall better results. 
 This is demonstrated, for example, in Figure \ref{fig:GSE-2SGR-LRT-casewise-curve}  which shows the average LRT distance behaviors for various casewise contamination values, $k$.

\afterpage{
\begin{table}[t!]
\centering
\footnotesize
\caption{Finite sample efficiency for random correlations. The sample size is $n=10p$. }\label{tab:GSE-2SGRE-efficiency}
\begin{tabular}{lccccccHHHHc}
  \hline
$p$ & MLE & Rocke & HSD & Snip & DMCDSc & UF- & UBF- & UF- & UBF- & DDC- & UBF-DDC- \\ 
 & & & & &  & GSE & GSE & GRE-C & GRE-C & GRE-C & GRE-C \\

  \hline
10 & 1.00 & 0.50 & 0.73 & 0.12 & 0.41 & 0.75 & 0.66 & 0.53 & 0.48 & 0.56 & 0.57 \\ 
  20 & 1.00 & 0.57 & 0.92 & 0.09 & 0.56 & 0.83 & 0.73 & 0.59 & 0.55 & 0.61 & 0.61 \\ 
  30 & 1.00 & 0.58 & 0.93 & 0.10 & 0.63 & 0.87 & 0.79 & 0.49 & 0.44 & 0.48 & 0.50 \\ 
  40 & 1.00 & 0.60 & 0.94 & 0.10 & 0.68 & 0.89 & 0.83 & 0.39 & 0.36  & 0.40 & 0.40 \\ 
  50 & 1.00 & 0.60 & 0.94 & 0.11 & 0.70 & 0.91 & 0.84 & 0.48 & 0.49 & 0.56 & 0.58 \\ 
   \hline
\end{tabular}
\end{table}
}

Table \ref{tab:GSE-2SGRE-efficiency} shows the finite sample relative efficiency under clean samples with random correlation for the considered robust estimates, taking the MLE average LRT distances as the baseline. The results for the AR$1(0.9)$ correlation are very similar and not shown here. As expected, UF-GSE show an increasing efficiency as $p$ increases while UBF-DDC-GRE-C have lower efficiency. Improvements can be achieved by using smaller $\alpha$ in the Rocke $\rho$ function with some trade-off in robustness. Results from this experiment are provided in the  supplementary material.

\afterpage{
\begin{table}[t!]
\centering
\footnotesize
\caption{Average ``CPU time" -- in seconds of a 2.8 GHz Intel Xeon -- evaluated using the \texttt{R} command, \texttt{system.time}. The sample size is $n=10p$.} \label{tab:GSE-2SGRE-computing-times}
\begin{tabular}{lcHHHc}
  \hline
$p$ & UF- & UBF- & UF- & UBF- & UBF-DDC- \\ 
& GSE & GSE & GRE-C & GRE-C & GRE-C \\
  \hline
10 & 0.7 & 1.1 & 0.1 & 0.2 & 0.2 \\ 
  20 & 7.7 & 11.0 & 1.2 & 1.7 & 1.7 \\ 
  30 & 34.5 & 45.6 & 5.4 & 6.3 & 6.4 \\ 
  40 & 120.5 & 144.9 & 14.5 & 16.9 & 17.1 \\ 
  50 & 278.4 & 338.0 & 33.0 & 37.0 & 37.8\\ 
   \hline
\end{tabular}
\end{table}
}

Finally, we compare the computing times of the two-step procedures. Table \ref{tab:GSE-2SGRE-computing-times} shows the average computing times over all contamination settings for various dimensions and for $n =10p$. The computing times for the two-step procedure have been substantially improved with the implementation of the faster initial estimator, EMVE-C.

\section{Real data example: small-cap stock returns data}\label{sec:GSE-example}

In this section, we consider the weekly returns from 01/08/2008 to 12/28/2010  for a portfolio of 20 small-cap stocks  from  \citet{martin:2013}. 

The  purpose of this example is  fourfold: first, to show that the classical MLE  and traditional robust procedures perform poorly on data affected by propagation of cellwise outliers; second, to show that the two-step procedures (e.g., UF-GSE) can provide better estimates by filtering large outliers; third, that the bivariate-filter version of the two-step procedure (e.g., UBF-GSE) provides even better estimates by flagging additional moderate cellwise outliers; and fourth,  that the two-step procedures that use GRE-C (e.g., 
UBF-GRE-C) can more effectively downweight some high-dimensional casewise outliers than those that use GSE (e.g., UBF-GSE), for this 20-dimensional dataset. Therefore, UBF-GRE-C provides the best results for this dataset. 

\begin{figure}[thb]
\centering 
\includegraphics[scale=0.44]{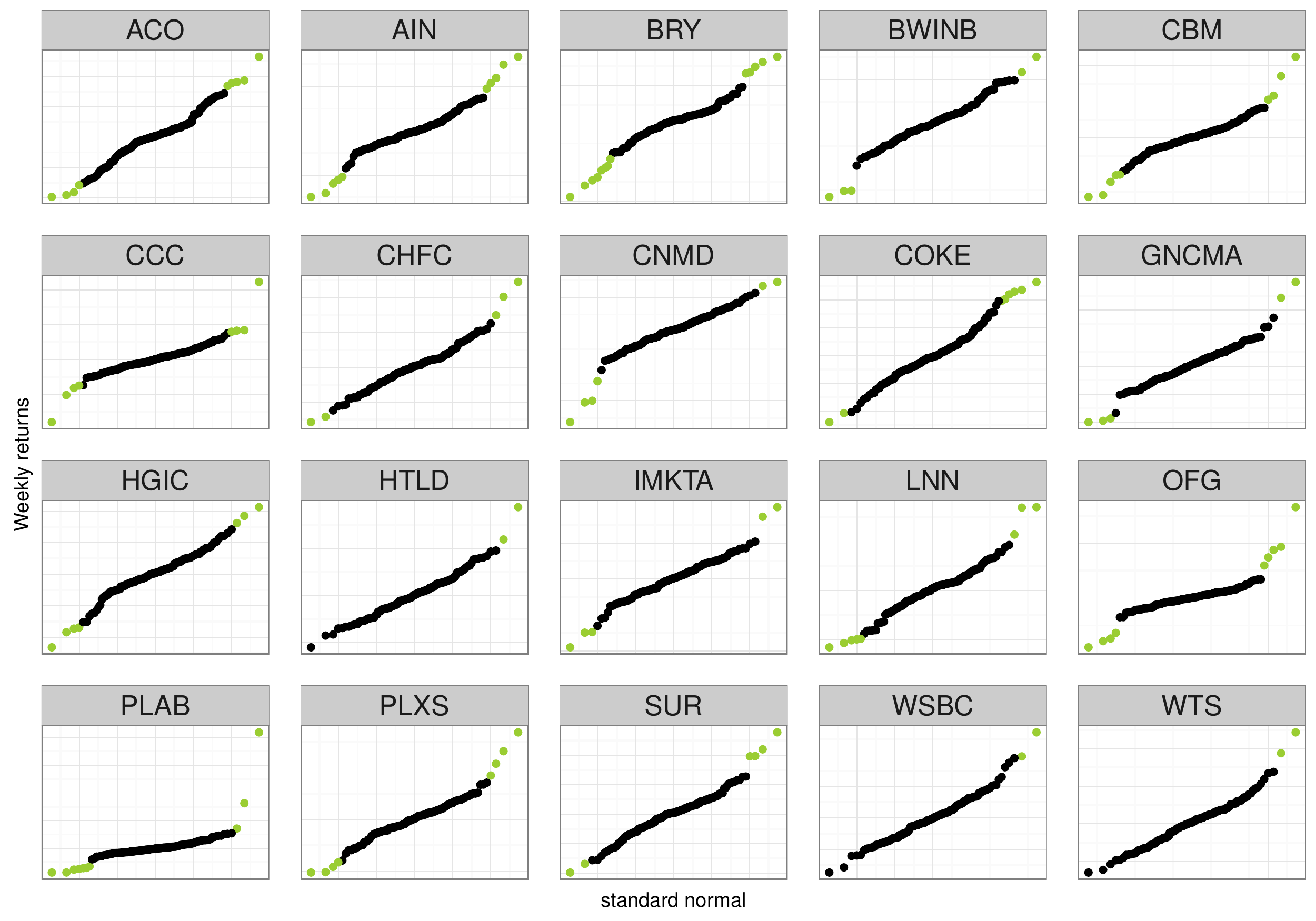} 
\caption{Normal quantile--quantile plots of weekly returns. Weekly returns that are three MAD's away from the coordinatewise-median are shown in green. }\label{fig:GSE-stock-qq}
\end{figure}

 Figure \ref{fig:GSE-stock-qq} shows the normal QQ-plots of the 20 small-cap stocks returns in the portfolio. The bulk of the returns in all stocks seem roughly normal, but large outliers are clearly present for most of these stocks. Stocks with returns lying more than three MAD's away from the coordinatewise-median (i.e., the large outliers) are shown in green in the figure. There is a total of 4.8\%  large cellwise outliers that propagate to 40.1\% of the cases.  Over 75\% of these weeks  correspond to the 2008 financial crisis.

\begin{figure}[t!]
\centering
\includegraphics[scale=0.65]{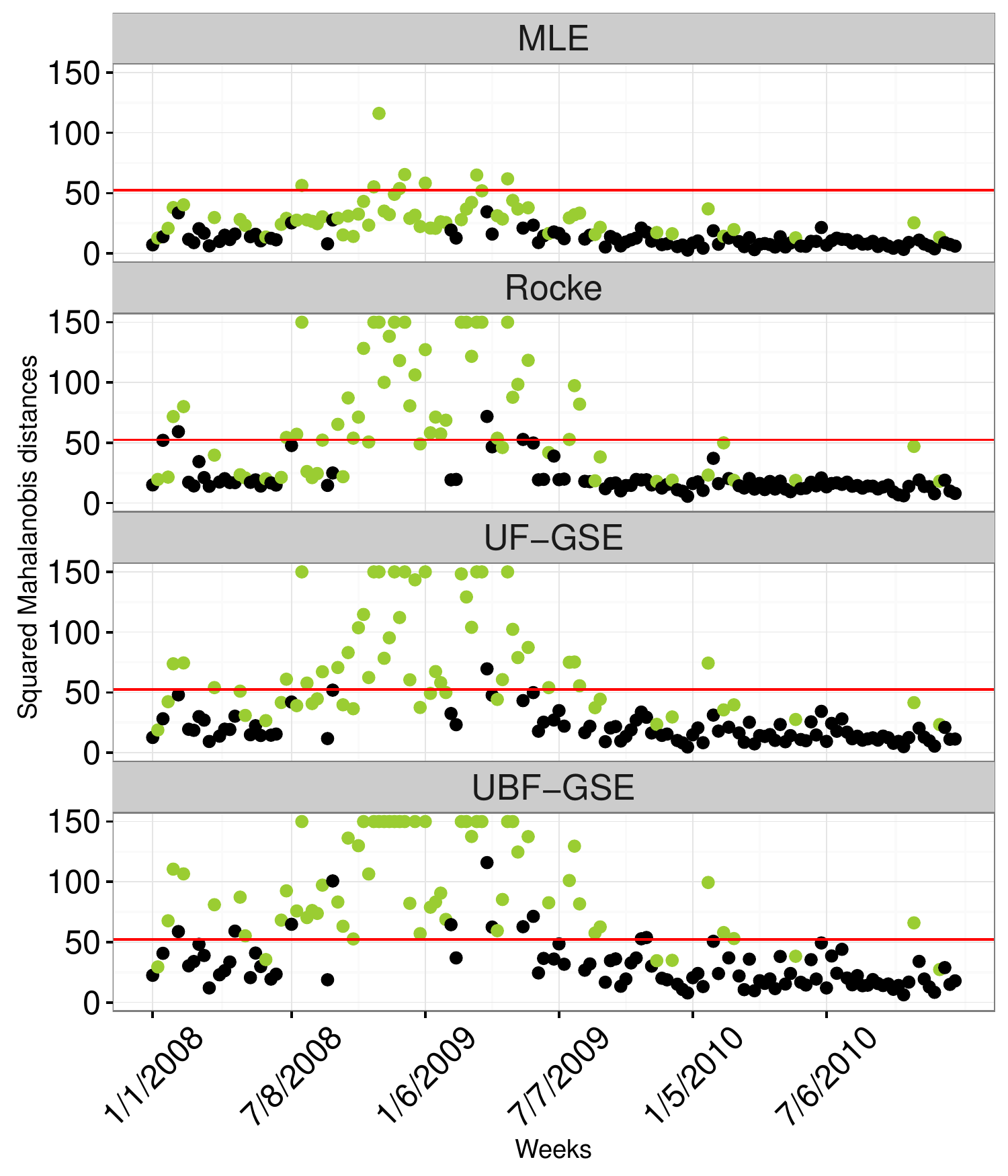}
\caption{
Squared Mahalanobis distances of the weekly observations in the small-cap asset returns data based on the MLE, the Rocke, the UF-GSE, and the UBF-GSE estimates. Weeks that contain one or more asset returns with values three MAD's away from the coordinatewise-median are in green. Large distances are truncated for better visualization.}\label{fig:GSE-stock-distances-Rocke-2SGS}
\end{figure}

\begin{figure}[t!]
\centering 
\includegraphics[scale=0.65]{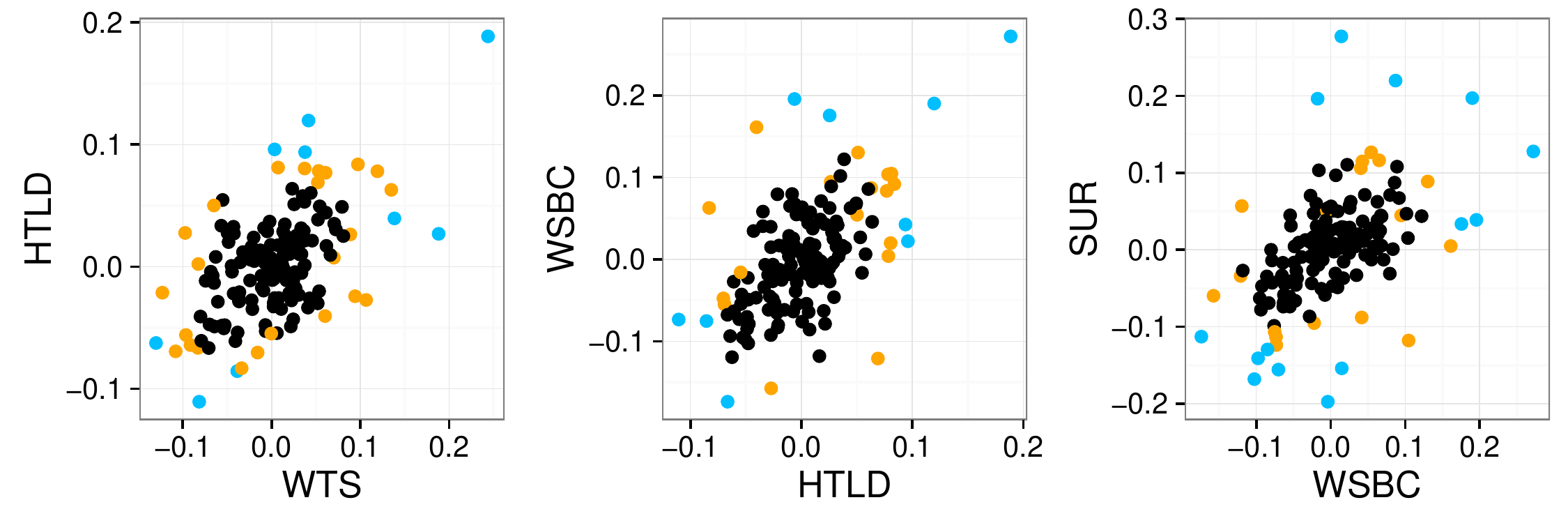} 
\caption{Pairwise scatterplots of the asset returns data for WTS versus HTLD, HTLD versus WSBC, and WSBC versus SUR.  Points with components  flagged by the univariate filter are  in blue. Points with components additionally flagged by the bivariate filter are in
orange.}\label{fig:GSE-stock-filtering-results}
\end{figure}

 Figure \ref{fig:GSE-stock-distances-Rocke-2SGS} shows the squared Mahalanobis distances of the 157 weekly observations based on four estimates: the MLE, the Rocke-S estimates, the UF-GSE, and the UBF-GSE.
 Weeks that contain  large cellwise outliers (asset returns with values three MAD's away from the coordinatewise-median) are in green.  From the figure, we see that the MLE and the Rocke-S estimates have failed to identify many of those weeks as MD outliers (i.e., failed to flag these weeks as having estimated full Mahalanobis distance exceeding the 99.99\% quantile chi-squared distribution with 20 degrees of freedom). The MLE misses all but seven of the 59 green cases.  The Rocke-S  estimate does slightly better but still misses one third of the  green cases. This is because it is severely affected by the large  cellwise outliers that propagate to  40.1\% of the cases.  The UF-GSE estimate also does a relatively poor job.  This may be due to the presence of several moderate cellwise outliers. In fact, 
Figure \ref{fig:GSE-stock-filtering-results} shows  the pairwise scatterplots for WTS versus HTLD, HTLD versus WSBC, and WSBC versus SUR with the results from the univariate and the bivariate filter. The points flagged by the univariate filter are in blue, and those flagged by the bivariate filter are in orange. We see that the bivariate filter has identified some additional cellwise outliers that are not-so-large marginally but become more visible when viewed together with other correlated components. These moderate cellwise outliers account for 6.9\% of the cells in the data and  propagate to 56.7\% of the cases. The final median weight assigned to these cases by UF-GSE and UBF-GSE are 0.50 and 0.65, respectively. By filtering the moderate cellwise outliers, UBF-GSE makes a more effective  use of the clean part of these  partly contaminated data points (i.e., the 56.7\% of the cases). As a result, UBF-GSE successfully flags all but five of the 59 green cases.

\begin{figure}[t!]
\centering
\includegraphics[scale=0.65]{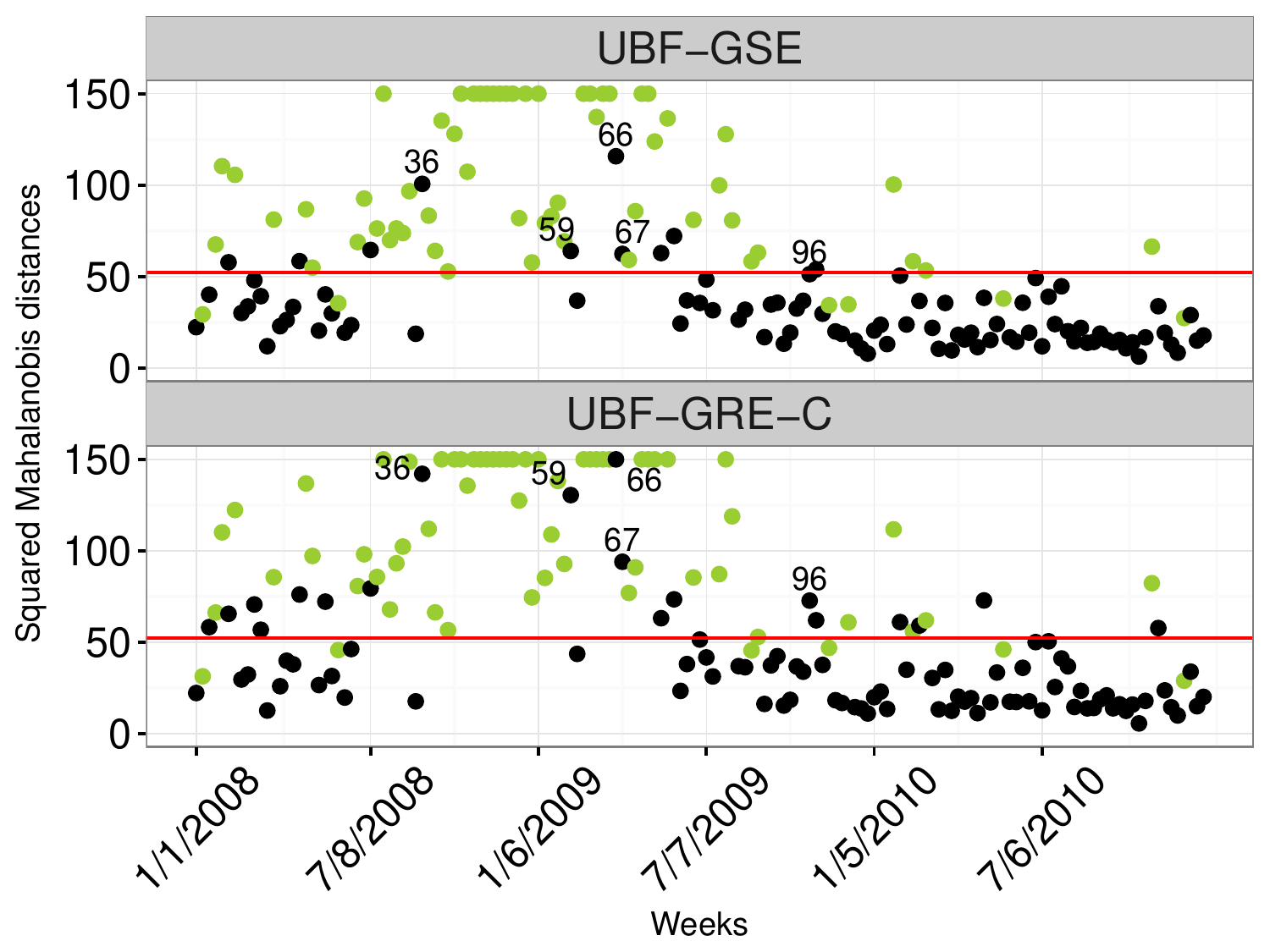}
\caption{
Squared Mahalanobis distances of the weekly observations in the small-cap asset returns data based on the UBF-GSE and the UBF-GRE-C estimates. Weeks that contain one or more asset returns with values three MAD's away from the coordinatewise-median are in green. }\label{fig:GSE-stock-distances-2SGR}
\end{figure}

Figure \ref{fig:GSE-stock-distances-2SGR} shows the squared Mahalanobis distances produced by UBF-GRE-C and  UBF-GSE, for comparison.  Here, we  see that UBF-GRE-C has missed only 3 of the 59 green cases, while UBF-GSE has missed 6 of the 59. UBF-GRE-C has also clearly flagged weeks 36, 59, and  66 (with final weights 0.6, 0.0, and 0.0, respectively) as casewise outliers.  In contrast,  UBF-GSE gives final weights 0.8, 0.5, and 0.5 to these cases. Consistent with our  simulation results, UBF-GSE has difficulty downweighting  some high-dimensional outlying cases on datasets of high dimension. 

In this example, UBF-GRE-C makes the most effective use of the clean part of the data and has the best outlier detecting performance among the considered estimates.

\section{Conclusions}\label{sec:GSE-conclusions}

In this paper, we overcome three serious limitations of UF-GSE. First, the estimator cannot deal with moderate cellwise outliers. Second, the estimator shows an incontrollable  increase in Gaussian efficiency, which is paid off by a serious decrease in robustness, for larger $p$. Third, the initial estimator (extended minimum volume ellipsoids, EMVE) used by GSE and UF-GSE  does not scale well in higher dimensions because it requires an impractically large number of subsamples to achieve a high breakdown point in  larger dimensions. 

To deal with also moderate cellwise outliers, we complement the univariate filter with a combination of bivariate filters (UBF-DDC). 
To achieve a controllable efficiency/robustness trade off  in higher dimensions, we replace the GSE in the second step with the Rocke-type GSE which we called it GRE. Finally,
to overcome the high computational cost of the EMVE, we introduce a clustering-based subsampling procedure. The proposed procedure is called UBF-DDC-GRE-C.

As shown by our simulation, UBF-DDC-GRE-C provides reliable results  for cellwise contamination when $\epsilon \le 0.05$ and $p \le 50$. For larger dimensions ($p > 50$), in our experience, the proposed estimator still performs well unless there is  a large fraction of small size cellwise outliers that evade the filter and propagate. 
Furthermore, UBF-DDC-GRE-C exhibits high robustness against moderate and large cellwise outliers, as well as casewise outliers in higher dimensions (e.g., $p > 10$). 
We also show via simulation studies that, in higher dimensions, estimators using the proposed subsampling with only $50$ subsamples can achieve equivalent performance than  the usual uniform subsampling with $500$ subsamples.

The proposed two-step procedure still has some limitation. 
As pointed out in the rejoinder in \citet{agostinelli:2014b}, the GSE in the second step does not handle well flat data sets, i.e., $n \approx 2p$. In fact, when $n \le 2p$, these estimators fail to exist (cannot be computed). This is also the case  for  GRE-C,  and for all the casewise robust estimators with breakdown point $1/2$. Our numerical experiments show that the proposed two-step procedure works well when $n \ge 5p$ but  not  as well when $2p < n < 5p$, depending on the amount of data filtered in the first step. In this situation, if much data are filtered leaving a small fraction of complete data cases, GSE and GRE may fail to converge \citep{danilov:2012, agostinelli:2014b}. This problem could be remedied by using graphical lasso \citep[GLASSO,][]{friedman:2008} to improve the conditioning of  the  estimates.

\appendix
\section{Proofs of Propositions}

\subsection{Proof of Proposition \ref{prop:2SGS-GY-asymptotic}}

The proof was available in \citet{agostinelli:2014}, but we provide a more detailed proof in the supplementary material for completeness.

\subsection{Proof of Proposition \ref{prop:2SGS-GY-asymptotic-bivariate}}\label{sec:2SGS-GY-asymptotic-proof-bivariate}

We need the following lemma for the proof.

\begin{lemma}\label{lem:2SGS-GY-asymptotic-bivariate}
Consider a sample of $p$-dimensional random vectors $\pmb X_1, \dots, \pmb X_n$. Also, consider a pair of multivariate location and scatter estimators $\pmb T_{0n}$ and $\pmb C_{0n}$. Suppose that $\pmb T_{0n}\to \pmb\mu_0$ and $\pmb C_{0n} \to \pmb\Sigma_0$ a.s.. Let $D_{i} = (\pmb X_i - \pmb T_{0n})^t \pmb C_{0n}^{-1}(\pmb X_i - \pmb T_{0n})$ and $D_{i} = (\pmb X_i - \pmb \mu_{0})^t \pmb \Sigma_{0}^{-1}(\pmb X_i - \pmb \mu_{0})$. Given $K < \infty$. 
For all $i=1,\dots, n$, if $D_{0i} \le K$ , then:
\begin{linenomath}
\[
D_{i} \to D_{0i} \quad \text{a.s..}
\]
\end{linenomath}
\end{lemma}

\begin{proof}[{\it Proof of Lemma \ref{lem:2SGS-GY-asymptotic-bivariate}}]
Note that 
\begin{linenomath}
\begin{align*}
|D_i - D_{0i}| &= |(\pmb X_i - \pmb T_{0n})^t \pmb C_{0n}^{-1} (\pmb X_i - \pmb T_{0n}) - ( \pmb X_i - \pmb\mu_0)^t \pmb\Sigma_0^{-1} (\pmb X_i - \pmb\mu_0) | \\
&= \begin{multlined}[t]
|((\pmb X_i - \pmb \mu_0) + (\pmb\mu_0 - \pmb T_{0n}))^t (\pmb\Sigma_0^{-1} + (\pmb C_{0n}^{-1}- \pmb\Sigma_0^{-1})) ((\pmb X_i - \pmb \mu_0 )+ (\pmb\mu_0 - \pmb T_{0n})) \\
- (\pmb X_i - \pmb\mu_0)^t \pmb\Sigma_0^{-1} (\pmb X_i - \pmb\mu_0)|
\end{multlined}
\\
&\le 
\begin{multlined}[t]
|(\pmb\mu_0 - \pmb T_{0n})^t\pmb\Sigma_0^{-1} (\pmb\mu_0 - \pmb T_{0n})|
+ |(\pmb\mu_0 - \pmb T_{0n})^t(\pmb C_{0n}^{-1}- \pmb\Sigma_0^{-1}) (\pmb\mu_0 - \pmb T_{0n})|\\
+ |2(\pmb X_i - \pmb \mu_0)^t \pmb\Sigma_0^{-1}(\pmb\mu_0 - \pmb T_{0n}) | + |2(\pmb X_i - \pmb \mu_0)^t (\pmb C_{0n}^{-1}- \pmb\Sigma_0^{-1})(\pmb\mu_0 - \pmb T_{0n}) | \\
+ |(\pmb X_i - \pmb \mu_0)^t(\pmb C_{0n}^{-1}- \pmb\Sigma_0^{-1}) (\pmb X_i - \pmb \mu_0)| 
\end{multlined}
\\
&= A_n + B_n + C_n + D_n + E_n.
\end{align*}
\end{linenomath}

By assumption, there exists $n_1$ such that for $n \ge n_1$ implies $A_n \le \varepsilon/5$ and $B_n \le \varepsilon/5$.

Next, note that
\begin{linenomath}
\[
\begin{multlined}[t]
|(\pmb X_i - \pmb\mu_0)^t\pmb\Sigma_0^{-1/2} \pmb y| = |\pmb y^t \pmb\Sigma_0^{-1/2} (\pmb X_i - \pmb\mu_0)| \\
\le || \pmb y || ||\pmb\Sigma_0^{-1/2} (\pmb X_i - \pmb\mu_0)|| 
= ||\pmb y|| \sqrt{(\pmb X_i - \pmb\mu_0)^t \pmb\Sigma_0^{-1}(\pmb X_i - \pmb\mu_0)} \le ||\pmb y|| \sqrt{K}.
\end{multlined}
\]
\end{linenomath}
So, there exists $n_2$ such that $n \ge n_2$ implies
\begin{linenomath}
\begin{align*}
C_n &= |2(\pmb X_i - \pmb \mu_0)^t \pmb\Sigma_0^{-1}(\pmb\mu_0 - \pmb T_{0n}) | \\
&= |2(\pmb X_i - \pmb \mu_0)^t \pmb\Sigma_0^{-1/2}\pmb\Sigma_0^{-1/2}(\pmb\mu_0 - \pmb T_{0n}) | \\
&\le 2 || \pmb\Sigma_0^{-1/2}(\pmb\mu_0 - \pmb T_{0n}) || \sqrt{K}\\
&\le \varepsilon/5.
\end{align*}
\end{linenomath}
Similarly, there exists $n_3$ such that $n \ge n_3$ implies
\begin{linenomath}
\begin{align*}
D_n &=|2(\pmb X_i - \pmb \mu_0)^t (\pmb C_{0n}^{-1}- \pmb\Sigma_0^{-1})(\pmb\mu_0 - \pmb T_{0n}) | \\
&= |2(\pmb X_i - \pmb \mu_0)^t \pmb\Sigma_0^{-1/2}\pmb\Sigma_0^{1/2}(\pmb C_{0n}^{-1}- \pmb\Sigma_0^{-1})(\pmb\mu_0 - \pmb T_{0n}) | \\
&\le 2 || \pmb\Sigma_0^{1/2}(\pmb C_{0n}^{-1}- \pmb\Sigma_0^{-1})(\pmb\mu_0 - \pmb T_{0n}) || \sqrt{K}\\
&\le \varepsilon/5.
\end{align*}
\end{linenomath}
Also, there exists $n_4$ such that $n \ge n_4$ implies
\begin{linenomath}
\begin{align*}
E_n &=|(\pmb X_i - \pmb \mu_0)^t(\pmb C_{0n}^{-1}- \pmb\Sigma_0^{-1}) (\pmb X_i - \pmb \mu_0)| \\
&=|(\pmb X_i - \pmb \mu_0)^t \pmb\Sigma_0^{-1/2}\pmb\Sigma_0^{1/2} (\pmb C_{0n}^{-1}- \pmb\Sigma_0^{-1}) (\pmb X_i - \pmb \mu_0)| \\
&\le ||\pmb\Sigma_0^{1/2} (\pmb C_{0n}^{-1}- \pmb\Sigma_0^{-1}) (\pmb X_i - \pmb \mu_0)|| \sqrt{K}\\
&\le||(\pmb C_{0n}^{-1}- \pmb\Sigma_0^{-1}) || \,\, ||\pmb\Sigma_0^{1/2} (\pmb X_i - \pmb \mu_0)|| \sqrt{K}\\
&\le||(\pmb C_{0n}^{-1}- \pmb\Sigma_0^{-1}) ||K\\
&\le \varepsilon/5.
\end{align*}
\end{linenomath}

Finally, let $n_5 = \max\{n_1, n_2, n_3, n_4\}$, then for all $i$,  $n \ge n_5$ implies
\begin{linenomath}
\[
|D_i - D_{0i}| \le \varepsilon/5 + \varepsilon/5 +\varepsilon/5 +\varepsilon/5+\varepsilon/5 = \varepsilon.
\]
\end{linenomath}
\end{proof}

\begin{proof}[{\it Proof of Proposition \ref{prop:2SGS-GY-asymptotic-bivariate}}]

Let $D_{0i} = (\pmb X_i - \pmb\mu_0)^t \pmb\Sigma_0^{-1} (\pmb X_i - \pmb\mu_0)$ and $D_{i} = (\pmb X_i - \pmb T_{0n})^t \pmb C_{0n}^{-1} (\pmb X_i - \pmb T_{0n})$. 
Denote  the empirical distributions of  $D_{01}, \dots, D_{0n}$ and $D_{1}, \dots, D_{n}$ by
\begin{linenomath}
\[
G_{0n}(t) = \frac{1}{n}\sum_{i=1}^n I\left( D_{0i} \le t \right) \quad \text{and} \quad
G_{n}(t) = \frac{1}{n}\sum_{i=1}^n I\left( D_{i} \le t \right).
\]
\end{linenomath}

Note that 
\begin{linenomath}
\begin{align*}
\left| G_n(t) - G_{0n}(t) \right| &= \left| \frac{1}{n}\sum_{i=1}^n I\left( D_{i} \le t \right) - \frac{1}{n}\sum_{i=1}^n I\left( D_{0i} \le t \right) \right| \\
&= 
\begin{multlined}[t]
\left| \frac{1}{n}\sum_{i=1}^n I\left( D_{i} \le t \right) I(D_{0i} > K) + \frac{1}{n}\sum_{i=1}^n I\left( D_{i} \le t \right) I(D_{0i} \le K) \right. \\
\left. - \frac{1}{n}\sum_{i=1}^n I\left( D_{0i} \le t \right) I( D_{0i} > K) - \frac{1}{n}\sum_{i=1}^n I\left( D_{0i} \le t \right) I( D_{0i} \le K) \right|
\end{multlined} \\
&\le \begin{multlined}[t]
\left| \frac{1}{n}\sum_{i=1}^n I\left( D_{i} \le t \right) I(D_{0i} > K) -  \frac{1}{n}\sum_{i=1}^n I\left( D_{0i} \le t \right) I( D_{0i}> K)   \right| \\
 + \left| \frac{1}{n}\sum_{i=1}^n I\left( D_{i} \le t \right) I(D_{0i} \le K)
- \frac{1}{n}\sum_{i=1}^n I\left( D_{0i} \le t \right) I( D_{0i} \le K) \right|
\end{multlined} \\
& = |A_n| + |B_n|.
\end{align*}
\end{linenomath}
We will show that $|A_n| \to 0$ and $|B_n| \to 0$ a.s..

Choose a large $K$ such that $P_{G_0}( D_{0} > K) \le \varepsilon/8$. By law of large numbers, there exists $n_1$ such that for $n \ge n_1$ implies $|\frac{1}{n}\sum_{i=1}^n   I( D_{0i} > K)  - P_{G_0}( D_{0} > K)| \le \varepsilon/8$ and 
\begin{linenomath}
\begin{align*}
|A_n| &=  \left| \frac{1}{n}\sum_{i=1}^n [I\left( D_{i} \le t \right)
-  I\left( D_{0i} \le t \right) ] I( D_{0i} > K)\right|  \\
&\le  \frac{1}{n}\sum_{i=1}^n | I\left( D_{i} \le t \right)
-  I\left( D_{0i} \le t \right) | I( D_{0i} > K)\\  
& \le  \frac{1}{n}\sum_{i=1}^n   I( D_{0i} > K)  \\
& \le P_{G_0}( D_{0} > K) + \varepsilon/8 \\
& \le \varepsilon/8  + \varepsilon/8 =  \varepsilon/4.
\end{align*}
\end{linenomath}

By assumption, we have from Lemma \ref{lem:2SGS-GY-asymptotic-bivariate} that  $D_i \to D_{0i}$ a.s. for all $i$ where $D_{0i} \le K$. Let $E_i = D_i - D_{0i}$. 
So, with probability 1, there exists $n_2$ such that $n \ge n_2$ implies that $-\delta \le  E_i  \le \delta$ for all $i$. 
Then, 
\begin{linenomath}
\begin{align*}
B_n &=  \frac{1}{n}\sum_{i=1}^n [ I\left( D_{i} \le t \right)  -  I\left( D_{0i} \le t \right)] I( D_{0i} \le K)    \\
&= \frac{1}{n}\sum_{i : D_{0i}\le K} [I\left( D_{i} \le t \right)  -  I\left( D_{0i} \le t \right) ]\\
&= \frac{1}{n}\sum_{i : D_{0i}\le K} [ I\left(  D_{0i} \le t - E_i \right)  -  I\left( D_{0i} \le t \right) ]\\
&\le \frac{1}{n}\sum_{i : D_{0i}\le K} [ I\left( D_{0i}  \le t + \delta \right)  -  I\left( D_{0i} \le t \right)  ]  \\
&\le \frac{1}{n}\sum_{i=1}^n [ I\left( D_{0i}  \le t + \delta \right)  -  I\left( D_{0i} \le t \right)].
\end{align*}
\end{linenomath}
Also, 
\begin{linenomath}
\begin{align*}
B_n &=  \frac{1}{n}\sum_{i : D_{0i}\le K} [ I\left(  D_{0i} \le t - E_i \right)  -  I\left( D_{0i} \le t \right) ]\\
&\ge \frac{1}{n}\sum_{i : D_{0i}\le K} [ I\left( D_{0i}  \le t - \delta \right)  -  I\left( D_{0i} \le t \right)  ]  \\
&\ge \frac{1}{n}\sum_{i=1}^n [ I\left( D_{0i}  \le t - \delta \right)  -  I\left( D_{0i} \le t \right)  ]  
\end{align*}
\end{linenomath}
Now, by the Gilvenko--Cantelli Theorem, with probability one there exists $n_3$ such that $n\ge n_3$ implies that $\sup_t |\frac{1}{n}\sum_{i=1}^n  I\left( D_{0i}  \le t + \delta \right) - G_0(t + \delta) | \le \varepsilon/16$, \\
$\sup_t |\frac{1}{n}\sum_{i=1}^n  I\left( D_{0i}  \le t - \delta \right) - G_0(t - \delta) | \le \varepsilon/16$,
 and $\sup_t |\frac{1}{n}\sum_{i=1}^n  I\left( D_{0i}  \le t \right) - G_0(t ) | \le \varepsilon/16$. 
Also, by the uniform continuity of $G_0$,  there
exists $\delta> 0$ such that $|G_0(t + \delta)-G_0(t)|\leq\varepsilon/8$ and $|G_0(t - \delta)-G_0(t)|\leq\varepsilon/8$.
Together,
\begin{linenomath}
\begin{alignat*}{3}
 \frac{1}{n}\sum_{i=1}^n  I\left( D_{0i} \le t - \delta \right)  -  I\left( D_{0i} \le t \right) &\le B_n  &&\le \frac{1}{n}\sum_{i=1}^n  I\left( D_{0i}  \le t + \delta \right)  -  I\left( D_{0i} \le t \right) \\
G_0(t - \delta) - \varepsilon/16  - G_0(t) - \varepsilon/16 &\le B_n &&\le  G_0(t + \delta) + \varepsilon/16  - G_0(t) + \varepsilon/16\\
(G_0(t - \delta) - G(t)) - \varepsilon/8 &\le B_n &&\le  ( G_0(t + \delta) -  G_0(t ) ) + \varepsilon/8 \\
 - \varepsilon/8 - \varepsilon/8 =-\varepsilon/4 &\le B_n &&\le  \varepsilon/8 + \varepsilon/8 = \varepsilon/4.
\end{alignat*}
\end{linenomath}

Finally, note that 
\begin{linenomath}
\[
G(t) - G_n(t) = ( G(t) - G_0(t)  )  + (G_0(t) - G_{0n}(t)) +  (G_{0n}(t) - G_n(t)).
\]
\end{linenomath}
Let $n_4 =\max\{n_1, n_2, n_3\}$, then $n \ge n_4$ implies
\begin{linenomath}
 \begin{align*}
\sup_{t > \eta} (G(t) - {G}_{n}(t))  &\le
\sup_{t > \eta} ( G(t) - G_0(t)  ) + \sup_{t > \eta} (G_0(t) - G_{0n}(t)) 
+ \sup_{t > \eta} (G_{0n}(t) - G_n(t)) \\
&\le (\varepsilon/4 + \varepsilon/4) + \varepsilon/16 + 0 \le \varepsilon.
\end{align*}
\end{linenomath}

\end{proof}

\section{Additional Tables from the Simulation Study in Section \ref{sec:GSE-MCresults}}\label{sec:appendix-filter-comparison}

\begin{table}[h!]
\centering
\footnotesize
\caption{Maximum average LRT distances under cellwise contamination. The sample size is $n=10p$. }
\begin{tabular}{llcHHHHHHHHccc}
  \hline
Corr. & $p$ & $\epsilon$ & MLE & Rocke & HSD & Snip & DMCDSc & UF- & UBF- & UF- & UBF- & DDC- & UBF-DDC- \\ 
& & & & & & &  & GSE & GSE & GRE-C & GRE-C  & GRE-C & GRE-C\\
\hline
 Random & 10 & 0 & 0.6 & 1.2 & 0.8 & 5.0 & 1.5 & 0.8 & 0.9 & 1.2 & 1.3 &  1.0 & 1.0\\
   &  & 0.02 & 114.8 & 1.2 & 2.3 & 6.9 & 1.6 & 1.2 & 1.4 & 1.3 & 1.4 & 1.1 & 1.1\\ 
   &  & 0.05 & 285.4 & 3.6 & 11.2 & 7.5 & 3.2 & 4.5 & 4.4 & 2.2 & 2.5 & 2.6 & 2.5\\ 
   & 20 & 0 & 1.1 & 2.0 & 1.2 & 11.5 & 2.0 & 1.3 & 1.5 & 1.9 & 2.0 & 1.8 & 1.8 \\ 
   &  & 0.02 & 146.1 & 2.7 & 10.6 & 13.9 & 2.6 & 4.0 & 4.4 & 2.9 & 3.0 & 2.5 & 2.5  \\ 
   &  & 0.05 & 375.9 & 187.2 & 57.1 & 15.5 & 9.3 & 11.0 & 11.1 & 8.0 & 8.2 & 7.7 & 7.3\\ 
   & 30 & 0 & 1.6 & 2.8 & 1.7 & 16.7 & 2.6 & 1.9 & 2.0 & 3.4 & 3.9 & 3.5 & 3.3\\ 
   &  & 0.02 & 179.0 & 23.1 & 22.6 & 18.5 & 4.4 & 5.8 & 6.3 & 5.4 & 5.9 & 5.3 & 5.0 \\ 
   &  & 0.05 & 475 & 380.5 & 123.1 & 20.8 & 13.7 & 14.2 & 14.8 & 12.3 & 13.4 & 14.2 & 13.3 \\  
   & 40 & 0 & 2.1 & 3.6 & 2.3 & 20.7 & 3.2 & 2.4 & 2.6 & 5.9 & 6.2 & 5.8 & 5.8 \\ 
   &  & 0.02 & 215.1 & 121.3 & 38.9 & 22.6 & 6.0 & 7.3 & 8.0 & 9.4 & 10.9 & 9.5 & 8.8 \\ 
   &  & 0.05 & $>$500 & $>$500 & 212.4 & 25.8 & 17.9 & 16.6 & 17.4 & 18.4 & 19.9 & 18.8 & 18.6 \\ 
   & 50 & 0 & 2.7 & 4.4 & 2.8 & 25.4 & 3.8 & 2.9 & 3.2 & 5.2 & 5.3  & 4.9 & 4.9 \\ 
   &  & 0.02 & 249.0 & 192.8 & 58.7 & 27.1 & 8.1 & 9.1 & 10.0 & 12.5 & 12.9 & 12.5 & 12.1 \\ 
   &  & 0.05 & $>$500 & $>$500 & 298.7 & 29.7 & 20.7 & 19.6 & 20.6 & 22.7 & 23.6 & 24.4 & 23.8 \\ 
   \hline
AR1(0.9) & 10 & 0 & 0.6 & 1.1 & 0.8 & 4.3 & 1.4 & 0.7 & 0.8 & 1.1 & 1.2 & 1.1 & 1.0 \\ 
   &  & 0.02 & 149.8 & 1.2 & 0.9 & 4.9 & 1.5 & 0.9 & 0.9 & 1.2 & 1.3 & 1.1 & 1.0 \\ 
   &  & 0.05 & 383.8 & 2.6 & 2.8 & 7.0 & 3.1 & 2.1 & 1.1 & 1.7 & 1.4 & 1.3 & 1.3 \\ 
   & 20 & 0 & 1.1 & 1.9 & 1.2 & 7.8 & 2.1 & 1.2 & 1.3 & 1.8 & 1.9 & 1.8 & 1.7 \\ 
   &  & 0.02 & 311.3 & 2.5 & 3.9 & 10.5 & 2.6 & 2.1 & 1.5 & 2.2 & 2.1 & 2.0 & 1.9 \\ 
   &  & 0.05 & $>$500 & $>$500 & 31.3 & 14.3 & 12.3 & 9.3 & 2.7 & 7.6 & 2.8 & 2.1 & 2.5 \\ 
   & 30 & 0 & 1.6 & 2.8 & 1.8 & 9.4 & 2.7 & 1.7 & 1.8 & 3.2 & 3.4 & 3.6 & 3.2 \\ 
   &  & 0.02 & 475.9 & 71.1 & 10.7 & 13.9 & 5.4 & 4.0 & 2.3 & 3.9 & 3.4 & 3.5 & 3.3 \\ 
   &  & 0.05 & $>$500 & $>$500 & 103.3 & 19.8 & 22.6 & 20.3 & 6.2 & 18.1 & 5.5 & 3.4 & 3.6 \\ 
   & 40 & 0 & 2.1 & 3.6 & 2.2 & 10.9 & 3.4 & 2.3 & 2.3 & 5.5 & 5.7 &  5.8 & 5.5 \\ 
   &  & 0.02 & $>$500 & 222.1 & 22.7 & 16.2 & 8.9 & 6.7 & 3.5 & 6.5 & 5.7 & 6.0 & 5.6 \\ 
   &  & 0.05 & $>$500 & $>$500 & 259.9 & 23.7 & 34.8 & 31.4 & 14.0 & 29.7 & 12.4 & 6.1 & 5.9 \\ 
   & 50 & 0 & 2.7 & 4.4 & 2.8 & 13.0 & 4.0 & 2.8 & 2.9 & 5.5 & 5.2 & 4.6 & 5.0 \\  
   &  & 0.02 & $>$500 & $>$500 & 43.3 & 18.9 & 12.8 & 9.7 & 4.9 & 9.7 & 6.4 & 6.4 & 7.8 \\ 
   &  & 0.05 & $>$500 & $>$500 & $>$500 & 28.9 & 46.5 & 42.8 & 22.6 & 40.8 & 20.4 & 7.9 & 8.9 \\ 
   \hline
\end{tabular}
\end{table}

\begin{table}[h!]
\centering
\footnotesize
\caption{Maximum average LRT distances under casewise contamination. The sample size is $n=10p$. }
\begin{tabular}{llcHHHHHHHHccc}
  \hline
Corr. & $p$ & $\epsilon$ & MLE & Rocke & HSD & Snip & DMCDSc & UF- & UBF- & UF- & UBF- & DDC- & UBF-DDC- \\ 
& & & & & & &  & GSE & GSE & GRE-C & GRE-C  & GRE-C & GRE-C\\
   \hline
 Random & 10 & 0 & 0.6 & 1.2 & 0.8 & 5.0 & 1.5 & 0.8 & 0.9 & 1.2 & 1.3 & 1.0 & 1.0 \\
   &  & 0.10 & 43.1 & 2.8 & 3.9 & 44.4 & 4.9 & 9.7 & 18.5 & 11.0 & 19.1 & 9.4 & 7.7  \\ 
   &  & 0.20 & 89.0 & 4.7 & 21.8 & 110.3 & 123.6 & 91.8 & 146.8 & 30.1 & 53.0 & 25.3 & 23.7 \\ 
   & 20 & 0 & 1.1 & 2.0 & 1.2 & 11.5 & 2.0 & 1.3 & 1.5 & 1.9 & 2.0 & 1.8 & 1.8 \\ 
   &  & 0.10 & 77.0 & 3.4 & 13.4 & 76.9 & 37.8 & 29.7 & 50.1 & 11.5 & 20.9 & 9.5 & 9.1  \\ 
   &  & 0.20 & 146.7 & 5.6 & 95.9 & 166.5 & 187.6 & 291.8 & 311.4 & 22.0 & 49.3 &  18.0 & 17.4 \\ 
   & 30 & 0 & 1.6 & 2.8 & 1.7 & 16.7 & 2.6 & 1.9 & 2.0 & 3.4 & 3.9 & 3.5 & 3.3 \\ 
   &  & 0.10 & 100.0 & 4.3 & 26.1 & 82.3 & 118.6 & 75.3 & 101.3 & 12.8 & 21.8 & 10.6 & 9.9 \\  
   &  & 0.20 & 200.7 & 7.4 & 297.7 & 220.9 & 268.4 & 415.5 & 445.2 & 21.7 & 47.6 &  18.7 & 16.9 \\ 
   & 40 & 0 & 2.1 & 3.6 & 2.3 & 20.7 & 3.2 & 2.4 & 2.6 & 5.9 & 6.2 & 5.8 & 5.8 \\ 
   &  & 0.10 & 125.9 & 5.2 & 46.3 & 101.6 & 130.6 & 140.2 & 168.8 & 18.6 & 29.5 & 17.7 & 16.2 \\ 
   &  & 0.20 & 252.4 & 9.1 & $>$500 & 186.2 & 340.1 & $>$500 & 579.9 & 22.7 & 52.3 & 21.2 & 19.5 \\ 
   & 50 & 0 & 2.7 & 4.4 & 2.8 & 25.4 & 3.8 & 2.9 & 3.2 & 5.2 & 5.3 & 4.9 & 4.9 \\ 
   &  & 0.10 & 150.3 & 5.9 & 80.0 & 121.9 & 139.5 & 258.1 & 228.8 & 27.5 & 43.4 & 21.2 & 17.6 \\ 
   &  & 0.20 & 303.1 & 10.0 & $>$500 & 224.3 & 407.7 & $>$500 & $>$500 & 24.2 & 64.8 & 23.7 & 23.0 \\ 
   \hline
AR1(0.9) & 10 & 0 & 0.6 & 1.1 & 0.8 & 4.3 & 1.4 & 0.7 & 0.8 & 1.1 & 1.2 & 1.1 & 1.0 \\ 
   &  & 0.10 & 43.1 & 2.8 & 1.7 & 20.2 & 2.9 & 3.7 & 4.3 & 3.1 & 3.6 & 3.0 & 2.9 \\ 
   &  & 0.20 & 88.9 & 4.8 & 8.7 & 49.7 & 29.7 & 50.8 & 50.1 & 7.2 & 8.4 & 6.8 & 6.9 \\ 
   & 20 & 0 & 1.1 & 1.9 & 1.2 & 7.8 & 2.1 & 1.2 & 1.3 & 1.8 & 1.9 & 1.8 & 1.7 \\ 
   &  & 0.10 & 77.0 & 2.8 & 4.7 & 43.8 & 14.8 & 12.9 & 14.9 & 3.5 & 4.3 & 3.3 & 3.3 \\ 
   &  & 0.20 & 146.6 & 5.3 & 35.3 & 113.0 & 87.6 & 260.5 & 193.9 & 7.3 & 10.5 & 6.0 & 6.0\\ 
   & 30 & 0 & 1.6 & 2.8 & 1.8 & 9.4 & 2.7 & 1.7 & 1.8 & 3.2 & 3.4 & 3.6 & 3.2 \\ 
   &  & 0.10 & 98.9 & 3.4 & 8.9 & 66.1 & 32.2 & 31.3 & 37.7 & 4.1 & 5.1 & 4.2 & 4.1 \\ 
   &  & 0.20 & 200.5 & 8.2 & 155.5 & 144.8 & 122.9 & 372.7 & 365.1 & 8.4 & 13.3 & 6.9 & 6.8 \\ 
   & 40 & 0 & 2.1 & 3.6 & 2.2 & 10.9 & 3.4 & 2.3 & 2.3 & 5.5 & 5.7  & 5.8 & 5.5 \\ 
   &  & 0.10 & 124.9 & 4.3 & 15.6 & 83.7 & 49.2 & 69.1 & 75.5 & 6.4 & 7.3 & 5.8 & 6.4 \\ 
   &  & 0.20 & 253.0 & 9.2 & 430.3 & 151.9 & 209.3 & 477.6 & 479.7 & 10.0 & 17.4 & 8.9 & 8.7 \\ 
   & 50 & 0 & 2.7 & 4.4 & 2.8 & 13.0 & 4.0 & 2.8 & 2.9 & 5.5 & 5.2 & 4.6 & 5.0 \\ 
   &  & 0.10 & 150.2 & 5.1 & 26.5 & 103.3 & 64.4 & 148.2 & 160.1 & 7.6 & 8.1 & 7.5 & 7.9 \\ 
   &  & 0.20 & 302.6 & 10.1 & $>$500 & 188.5 & 276.0 & $>$500 & $>$500 & 11.0 & 21.2 & 10.0 & 8.8 \\ 
   \hline
\end{tabular}
\end{table}

\clearpage
\newpage
\section{Supplementary Materials}

Additional simulation results and related supplementary material referenced in the article can be found in a separate document, ``Supplementary Material".


\begin{footnotesize}
\section*{Acknowledgement}
\noindent Ruben Zamar and Andy Leung research were partially funded by the Natural Science and Engineering Research Council of Canada.
\end{footnotesize}

\bibliographystyle{elsarticle-harv} 
\bibliography{UBF-GRE}

\end{document}


\title{ {\bf Supplementary Material for: Multivariate Location and Scatter Matrix Estimation Under Cellwise and Casewise Contamination } }
\author[1]{{\bf Andy Leung}}
\author[2]{{\bf Victor J. Yohai}}
\author[1]{{\bf Ruben H. Zamar}}
\affil[1]{{\small Department of Statistics, University of British Columbia, 3182-2207 Main Mall, Vancouver, British Columbia V6T 1Z4, Canada}}
\affil[2]{{\small Departamento de Matem\'atica, Facultad de Ciencias Exactas y Naturales, Universidad de Buenos Aires, Ciudad Universitaria, Pabell\'on 1, 1426, Buenos Aires, Argentina}}
\renewcommand\Authands{, }

\maketitle

\begin{abstract}
This supplementary material contains all the proofs, additional simulation results, and related supplementary material referenced in the article ``Robust Estimation of Multivariate Location and Scatter in the Presence of Cellwise and Casewise Contamination". 
\end{abstract}

\section{Proof of Proposition 1}

The proof was available in \citet{agostinelli:2014}, but we provide a more detailed proof here for completeness.

\begin{proof}
Without loss of generality, set $\mu_0 = 0$ and $\sigma_0 = 1$. 
Let $Z_{0i} = \frac{X_i - \mu_0}{\sigma_0} = X_i$ and $Z_{i} =  \frac{X_i - T_{0n}}{S_{0n}}$.
Denote  the empirical distributions of  $Z_{01}, \dots, Z_{0n}$ and $Z_{1}, \dots, Z_{n}$ by
\begin{linenomath}
\[
F_{0n}^+(t) = \frac{1}{n}\sum_{i=1}^n I\left( |Z_{0i}| \le t \right) \quad \text{and} \quad
F_{n}^+(t) = \frac{1}{n}\sum_{i=1}^n I\left( |Z_{i}| \le t \right).
\]
\end{linenomath}

By assumption,  with probability one, there exists $n_1$ such that $n \ge n_1$ implies $0 < 1 - \delta \le S_{0n} \le 1 + \delta$ and $-\delta \le T_{0n} \le \delta$, and we have 
\begin{linenomath}
\[
\begin{aligned}
F_{n}^+(t) &= \frac{1}{n}\sum_{i=1}^n I\left(- t \le Z_{i} \le t  \right) =  \frac{1}{n}\sum_{i=1}^n I\left(- t \le \frac{X_i - T_{0n}}{S_{0n}} \le t  \right) \\
&= \frac{1}{n}\sum_{i=1}^n I\left(- t S_{0n} + T_{0n}  \le X_{i} \le t S_{0n} + T_{0n} \right)  \\
&\ge \frac{1}{n}\sum_{i=1}^n I\left(- t (1 - \delta) + T_{0n}  \le X_{i} \le t (1 - \delta) + T_{0n} \right) \\
&\ge \frac{1}{n}\sum_{i=1}^n I\left(- t (1 - \delta) + \delta  \le X_{i} \le t (1 - \delta) - \delta \right) \\
&= \frac{1}{n}\sum_{i=1}^n I\left( |X_{i}| \le t(1 - \delta) - \delta \right) = F_{0n}^+(t(1 - \delta) - \delta).
\end{aligned}
\]
\end{linenomath} 
Now, by the Glivenko--Cantelli Theorem, with probability one there
exists $n_2$ such that $n\geq n_2$ implies that $\sup_{t}| F_{0n}^+(t) -F_0^+(t)|\leq\varepsilon/2$. Also, by the uniform continuity of $F_0^+$, given $\varepsilon>0,$ there
exists $\delta> 0$ such that $|F_0^+(t(1-\delta) - \delta)-F_0^+(t)|\leq\varepsilon/2$.

Finally, note that 
\begin{linenomath}
\begin{align*}
{F}_{n}^+(t)  &  \geq  F_{0n}^+(t(1-\delta) -\delta)\\
&  =\left(  F_{0n}^+(t(1-\delta) -\delta)-F_0^+(t(1-\delta
)-\delta)\right) \\
&  \qquad+(F_0^+(t(1-\delta) -\delta)-F_0^{+}(t))+(F_0^{+}(t)- F^+(t))+F^{+}(t).
\end{align*}
\end{linenomath}
Let $n_3=\max(n_1,n_2)$, then $n\geq n_3$ imply
\begin{linenomath}
\begin{align*}
\sup_{t> \eta}( F^+(t)-{F}_{n}^+(t))  &
\leq\sup_{t> \eta}\left| F_0^+(t(1-\delta) -\delta)-  F_{0n}^{+}(t(1-\delta)
-\delta)\right| \\
& \qquad+\sup_{t> \eta}\left| F_0^{+}(t) - F_0^{+}(t(1-\delta)-\delta) \right| 
 +\sup_{t> \eta}(F^{+}(t)-F_0^{+}(t))\\
&  \leq \frac{\varepsilon}{2} + \frac{\varepsilon}{2} + 0 = \varepsilon.
\end{align*}
\end{linenomath}
This implies that $n_{0}/n\rightarrow0$ a.s..
\end{proof}

\section{Performance comparison between GSE and GRE}

We conduct a simulation study to compare the standalone performances of the second steps (i.e., the estimation step) in the two-step S-estimators: GRE-C starting from EMVE-C  versus GSE starting from EMVE. 
%

We consider clean and casewise contaminated samples from a $N_p(\pmb{\mu_0},\pmb{\Sigma_0})$ distribution with $p=10,20,\dots,50$ and $n=10p$. The simulation mechanisms are the same as that of Section 5, but in addition, 
 5\% of the cells in the generated samples are randomly selected and assigned a missing value. The number of replicates is $N=500$.

\begin{figure}[t!]
\centering
\includegraphics[scale=0.5]{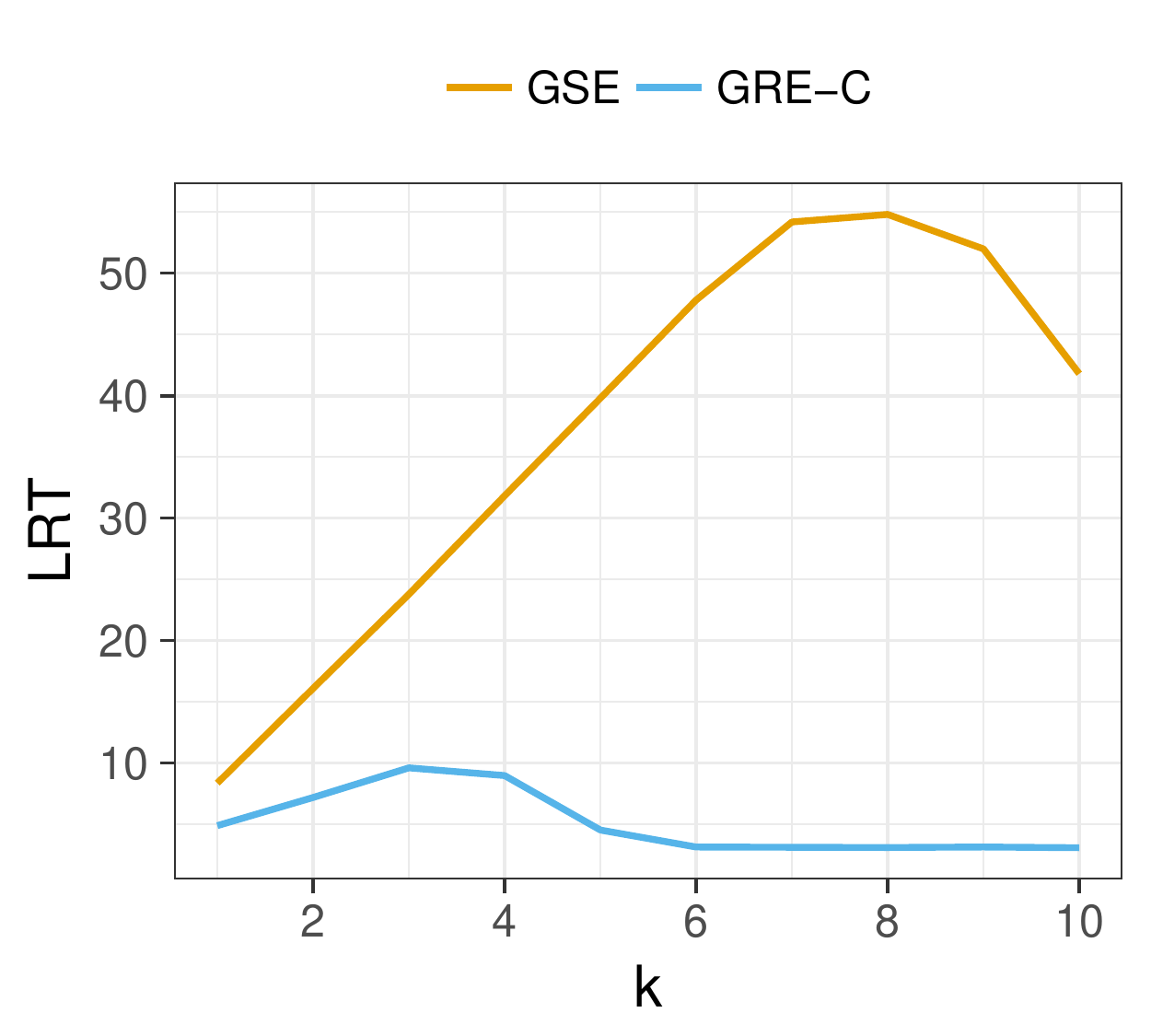}
\caption{Average LRT distances for various contamination sizes, $k$, for random correlations under 10\% casewise contamination. The dimension is $p=30$ and the sample size is $n=10p$.}\label{fig:GSE-LRT-curve-GRE}
\end{figure}

Table \ref{tab:GSE-robustness-all} shows the maximum average LRT distances from the true correlation matrices among the considered contamination sizes and, for brevity, shows only the values for random correlations. EMVE is capable of dealing small fraction of outliers with 500 subsamples, but breaks down when the fraction gets larger, and brings down the performance of GSE. EMVE-C with more refined subsampling procedure and larger subsample sizes shows better performance than EMVE, even for relatively larger fraction of outliers. Overall, GRE performs better than GSE. The Rocke $\rho$ function used in GRE is capable of giving smaller weights to points that are moderate-to-large distances from the main mass of points \citep{rocke:1996}; see, for example, 
Figure \ref{fig:GSE-LRT-curve-GRE} that shows the average LRT distance behaviors for $10\%$ contamination for dimension $30$ and sample size $300$ data. In the figure,  we see that GRE outperforms GSE for moderate sizes contamination points, as expected.

\begin{table}[t!]
\centering
\footnotesize
\caption{Maximum average LRT distances. The sample size is $n=10p$.}\label{tab:GSE-robustness-all}
\begin{tabular}{llcccc}
  \hline
$p$ & $\epsilon$ & EMVE & GSE & EMVE-C   & GRE-C \\ 
  \hline
10 & 0.10 & 8.7 & 4.6 & 17.3 & 10.9  \\ 
   & 0.20 & 81.4 & 84.8 & 43.4 & 36.1  \\ 
  20 & 0.10 & 20.8 & 24.1 & 9.2 & 8.1 \\ 
   & 0.20 & 123.0 & 156.8 & 13.1 & 14.9  \\ 
  30 & 0.10 & 31.2 & 54.8 & 13.4 & 9.4  \\ 
   & 0.20 & 299.1 & 223.2 & 24.3 & 16.0  \\ 
  40 & 0.10 & 77.5 & 80.7 & 21.9 & 12.2 \\ 
   & 0.20 & 511.8 & 287.9 & 43.2 & 17.1 \\ 
  50 & 0.10 & 172.5 & 125.1 & 29.4 & 16.5  \\ 
   & 0.20 & 644.3 & 349.8 & 60.2 & 26.3  \\ 
   \hline
\end{tabular}
\end{table}

Table \ref{tab:GSE-efficiency} shows the finite sample relative efficiency under clean samples, taking the classical EM estimator as the baseline. As expected, GSE shows an increasing efficiency as $p$ increases. GRE, overall, has lower efficiency.

\afterpage{
\begin{table}[t!]
\centering
\footnotesize
\caption{Finite sample efficiency. The sample size is $n=10p$.}\label{tab:GSE-efficiency}
\begin{tabular}{lcccc}
  \hline
$p$ & EMVE & GSE &  EMVE-C & GRE-C \\ 
  \hline
10 & 0.24 & 0.89 & 0.26 & 0.54  \\ 
  20 & 0.30 & 0.95 & 0.30 & 0.59 \\ 
  30 & 0.34 & 0.98 & 0.33 & 0.58  \\ 
  40 & 0.35 & 0.98 & 0.34 & 0.47 \\ 
  50 & 0.37 & 0.99 & 0.35 & 0.48 \\ 
     \hline
\end{tabular}
\end{table}
}

\section{Efficiency of GRE and tuning parameter $\alpha$ in Rocke-$\rho$ function}

The tuning parameter $\alpha$ in the Rocke-$\rho$ function in $\gamma$ in 
(10) is chosen small to control the efficiency. In this chapter, we used the conventional choice $\alpha = 0.05$, as seen to achieve reasonable efficiency while achieving high robustness. Here, we explore the performance of GRE-C with smaller values of $\alpha$. We repeat the simulation study as in Section 5 for $p=10,30,50$ and $n=10p$. The number of replicates is $N=30$. Table \ref{tab:GSE-GRE-alpha-efficiency} reports the finite sample efficiency and maximum average LRT distances under 20\% casewise contamination. In general, higher efficiency can be achieved using smaller values of $\alpha$, but with the cost of some loss in robustness.

\begin{table}[h!]
\centering
\footnotesize
\caption{Finite sample efficiency and maximum average LRT distances for GRE-C with various values of $\alpha$. The sample size is $n=10p$.}\label{tab:GSE-GRE-alpha-efficiency}
\begin{tabular}{lcccccc}
\hline
$p$ & \multicolumn{3}{c}{Efficiency, clean data} & \multicolumn{3}{c}{Max LRT, 20\%  casewise} \\ \cmidrule(lr){2-4} \cmidrule(lr){5-7}
  & $\alpha=0.05$ & $\alpha=0.01$ & $\alpha=0.001$ & $\alpha=0.05$ & $\alpha=0.01$ & $\alpha=0.001$ \\ 
\cmidrule(lr){1-4}\cmidrule(lr){5-7}
10 & 0.54 & 0.67 & 0.67 & 33.1 & 32.1 & 32.1 \\ 
30 & 0.58 & 0.85 & 0.95 & 16.0 & 20.2 & 28.7 \\ 
50 & 0.55 & 0.58 & 0.93 & 27.1 & 28.1 & 47.7 \\
\hline
\end{tabular}
\end{table}

\bibliographystyle{elsarticle-harv} 
\bibliography{UBF-GRE}